# ASYMPTOTIC GENEALOGY OF A CRITICAL BRANCHING PROCESS

BY LEA POPOVIC

*University of California, Berkeley*

Consider a continuous-time binary branching process conditioned to have population size $n$ at some time $t$, and with a chance $p$ for recording each extinct individual in the process. Within the family tree of this process, we consider the smallest subtree containing the genealogy of the extant individuals together with the genealogy of the recorded extinct individuals. We introduce a novel representation of such subtrees in terms of a point-process, and provide asymptotic results on the distribution of this point-process as the number of extant individuals increases. We motivate the study within the scope of a coherent analysis for an a priori model for macroevolution.

**1. Introduction.** The use of stochastic models in the theory of macroevolution (origin and extinction of species) has been common practice for many years now. Stochastic models have been used to recreate phylogenetic trees of extant taxa from molecular data, and to recreate the time series of the past number of taxa from the fossil record. However, only a few attempts have been made to make the two analyses consistent with each other. Instead of studying data-motivated models (which are scientifically more realistic for specific applications), the purpose of this paper is to study a purely random model that can accommodate such a coherent analysis. We study a mathematically fundamental stochastic model which allows for inclusion of both the extant and the fossil types of data in one analysis.

A significant interest in evolutionary biology is devoted to reconstructing phylogenies based on available data on the extant species (of molecular or other type). The assumption is that in the distant past there was a common ancestor from which the extant species evolved according to some (stochastic) evolutionary model. One then tries to find the ancestral history (genealogy) of the extant species which optimizes some "best fit" criterion. Using shapes of such phylogenetic trees one then hopes to make some inference on the diversification properties of the evolutionary process. We stress here that the number of extant species is given, although the time from the origin of the evolutionary process to the present often is inferred from the data as well. On the other hand, inference of diversification rates based on fossil data mostly makes use of time series analyses of fossil counts. Once again one assumes an underlying (stochastic) evolutionary process, then tries to use fluctuations in the time series of fossil counts to make inferences on the diversification rates of the process. For the most part, however, only a small fraction of the species are retained within the fossil record, with variation in sampling rates over time. It is often hard to estimate precisely the proportion of species retained as fossils, although its variability may be reasonably captured by considering the sampling rate to be random as well.

The motivation for this paper was to consider a stochastic process which would incorporate both sets of data within one evolutionary model, and to present results describing the genealogy of the fossil record as well. Our basic premises are the given information on the number of extant species, the amount of time from the origin of the process to the present and the chance of species to be retained in the fossil record. It is subsequently possible to randomize the amount of time from the origin to the present day, as well as the chance of being retained in the fossil record. Details of such randomization under a reasonable choice of priors can be found in [17].

The model we propose is the continuous-time critical branching process. The reasons for our choice are the following. If one is to consider a model in which extinctions and speciations are random without systematic tendencies for the number of species to increase or decrease, then for a branching process this translates into the criticality of the process (the average number of offspring of each individual is 1). Such a model corresponds to one general view in evolutionary biology that (except for mass extinctions and their aftermath) the overall number of species does not have exponential growth nor an exponential decrease.

The fundamental critical branching processes previously employed in evolutionary models have drawbacks that exclude their use in our proposed study. The basic evolution model is the Yule process [20], the elementary continuous-time pure birth process [the process starts with one individual, each individual gives birth to offspring according to a Poisson(rate 1) process]. One can clearly not employ this model, as it a priori does not involve the extinction of species, hence does not allow for inclusion of the fossil record. The next candidate model, which includes the extinction of individuals, is the basic neutral model used in population genetics. The Moran model [7] is the process of uniformly random speciations and extinctions of individuals in a population of a fixed size [the total number of individuals is a fixed number, each individual lives for an Exponential(mean 1) lifetime, at the end of which it is replaced by an offspring chosen uniformly at random from the total population including itself]. One can consider this process as having persisted from a distant past to the present, giving implicitly a genealogical tree of the extant individuals. Asymptotically in the total population size (with suitable rescaling) this genealogical process (backward in time) is Kingman's coalescent model. Although it is possible to make modifications of this model to allow for nonconstant population size [11], this unfortunately requires an a priori assumption on the evolution of the total population size in time.

We are interested in considering a group of species that have some common ancestor at their origin. This corresponds to the practice in evolutionary biology of considering monophyletic groups. In this sense, the critical continuous-time binary branching process, in which individuals live for an Exponential(mean 1) time during which they produce offspring at Poisson(rate 1) times, is the natural basic model for the given purpose. We want to study the genealogical structure of the process conditioned on its population size at a given time $t$. By genealogical structure we mean a particular subtree of the branching process family tree. We consider all the extant individuals at time $t$, and the subset of the extinct individuals each having independent chance $p$ of being sampled into the record. The genealogical subtree is the smallest one containing all the common ancestors of the extant individuals and all the sampled extinct individuals. We introduce a point-process representation of this genealogical subtree, with a convenient graphical interpretation, and derive its law. Our main result is the asymptotic behavior of such point-processes and their connection to a conditioned Brownian excursion.

The relationship between random trees and Brownian excursions has been much explored in the literature. We note only a small selection that is directly relevant to the work in this paper. Neveu and Pitman [14, 15] and Le Gall [12] noted the appearance of continuous-time critical branching processes embedded in the structure of a Brownian excursion. Abraham [1] and Le Gall [13] considered the construction of an infinite tree within a Brownian excursion, which is in some sense a limit of the trees from the work of Neveu and Pitman. The convergence of critical branching processes conditioned on total population size to a canonical tree within a Brownian excursion (the *continuum random tree*) was introduced by Aldous [3]. We state a connection of the asymptotic results in this paper with the above-mentioned results.

Some aspects of the genealogy of critical Galton–Watson trees conditioned on nonextinction have been studied by Durrett [6], without the use of random trees. It has also been studied within the context of superprocesses (see, e.g., [13]). We further note that, in the biological literature, models of evolution have been made on each level of taxonomy (species, genera etc.) separately, while it is certainly desirable to insure hierarchical consistency between them. A natural way to extend our analysis to the next taxonomic level is to superimpose on the branching process a random process of marks distinguishing some species as originators of a new higher taxon. In collaboration with Aldous, we have pursued this study in [4], as part of a larger project on coherent and consistent stochastic models for macroevolution.

As a last remark, we note that, as implied by general convergence results on critical branching processes ([3] and many others), the same asymptotic genealogical process obtained here should invariably hold in general for any critical branching process with finite offspring variance.

The paper is structured in two parts. In Section 2 we give a precise definition of the genealogical point-process representing the common ancestry of the extant individuals and provide its exact law and asymptotic behavior (Theorem 5). Then, in Section 3 we give the definition of the corresponding genealogical point-process that includes the sampled extinct individuals as well, and we provide its exact law and asymptotic behavior (Theorem 9).

**2. Genealogy of extant individuals.** Let $\mathcal{T}$ be a continuous-time critical branching process, with initial population size 1. In such a process each individual has an Exponential(rate 1) lifetime, in the course of which it gives birth to new individuals at Poisson(rate 1) times, with all the individuals living and reproducing independently of each other. Let $\mathcal{T}_{t,n}$ be the process $\mathcal{T}$ conditioned to have population size $n$ at time $t$. We shall use the same notation ($\mathcal{T}$ and $\mathcal{T}_{t,n}$) for the random trees with edge-lengths that are the family trees of these processes.

We depict these family trees as rooted planar trees with the following conventions. Each individual is represented with a set of edges whose total length is equal to that individual's lifetime. Each birth time of an offspring corresponds to a branch-point in the parent's edge, with the total length of the parent's edge until the branch-point equal to the parent's age at this time. The new individual is then represented by the edge on the right, while the parent continues in the edge on the left. Such trees are identified by their shape and by the collection of the birth times and lifetimes of individuals. We shall label the vertices in the tree in a depth-first search manner. An example of a random tree realization of $\mathcal{T}_{t,n}$ is shown in Figure 1(a).

REMARK 1. The random tree $\mathcal{T}$ we defined is almost the same as the family tree of a continuous-time critical binary-branching Galton–Watson process. The difference between the two is only in the identities of the individuals. If, in the Galton–Watson process, at each branching event with two offspring we were to impose the identification of the left offspring with its parent, the resulting random tree would be the same as the family tree of our branching process $\mathcal{T}$.

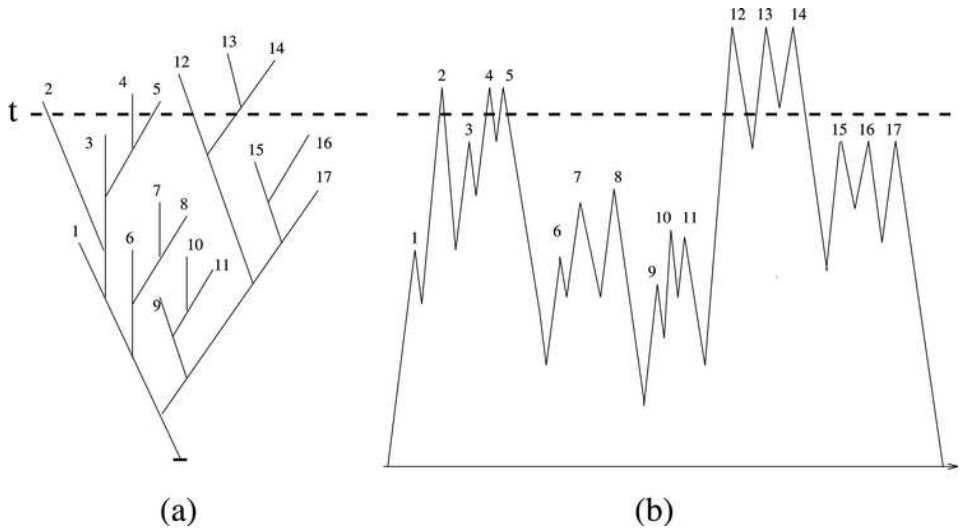

FIG. 1. (a) *A realization of the tree $\mathcal{T}_{t,n}$ whose population at time $t$ is $n = 5$; the leaves are labeled in depth-first search manner.* (b) *The contour $\mathcal{C}_{\mathcal{T}_{t,n}}$ process of the tree $\mathcal{T}_{t,n}$; each local maximum of $\mathcal{C}_{\mathcal{T}_{t,n}}$ corresponds to the height of a leaf of $\mathcal{T}_{t,n}$.*

Let $\mathcal{C}_\mathcal{T}$ be the contour process induced by the random tree $\mathcal{T}$. The contour process of a rooted planar tree is a continuous function giving the distance from the root of a unit-speed depth-first search of the tree. Such a process starts at the root of the tree, traverses each edge of the tree once upward and once downward following the depth-first search order of the vertices and ends back at the root of the tree. The contour process consists of line segments of slope $+1$ (the rises) and line segments of slope $-1$ (the falls). The unit speed of the traversal insures that the height levels in the process are equivalent to distances from the root in the tree, in other words to the times in the branching process. The contour process induced by the random tree $\mathcal{T}_{t,n}$ depicted in Figure 1(a) is $\mathcal{C}_{\mathcal{T}_{t,n}}$ shown in Figure 1(b). For a formal definition of planar trees with edge lengths, contour processes and their many useful properties one can consult the recent lecture notes of Pitman ([16], Section 6.1).

Let the *genealogy* of extant individuals at time $t$ be defined as the smallest subtree of the family tree which contains all the edges representing the ancestry of the extant individuals. The genealogy of extant individuals at $t$ in $\mathcal{T}_{t,n}$ is thus an $n$-leaf tree, which we denote by $\mathcal{G}(\mathcal{T}_{t,n})$. Figure 2(a) shows the genealogical subtree of the tree from Figure 1(a). We now introduce a novel point-process representation of this genealogical tree $\mathcal{G}(\mathcal{T}_{t,n})$. Thus we get an object that is much simpler to analyze and gives much clearer asymptotic results than if made in the original space of trees with edge-lengths.

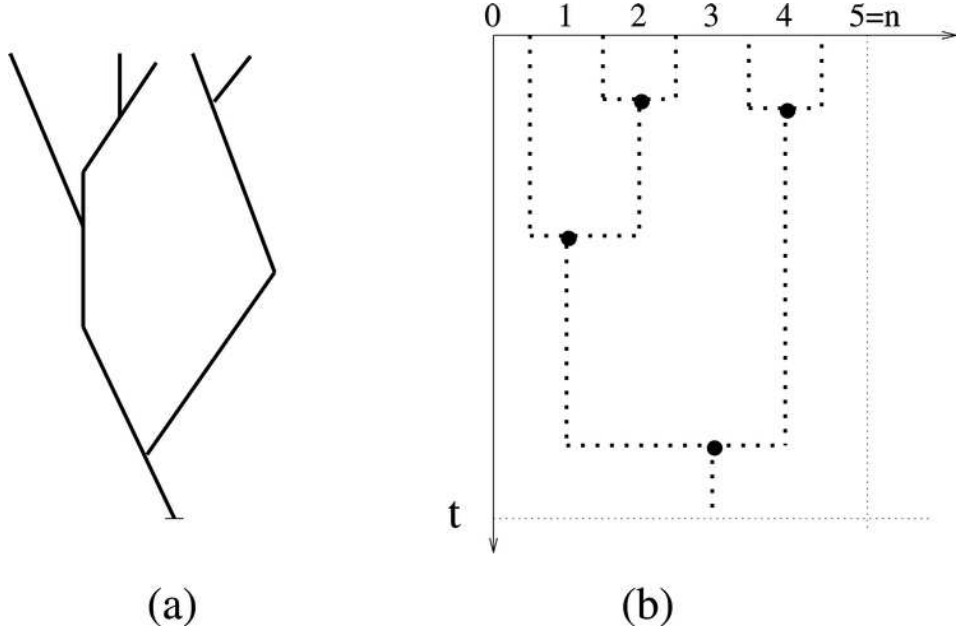

FIG. 2. (a) *The genealogical tree $\mathcal{G}(\mathcal{T}_{t,n})$ of the extant individuals at time $t$.* (b) *The point-process $\Pi_{t,n}$ representation of $\mathcal{G}(\mathcal{T}_{t,n})$* [*the dotted lines show the simple reconstruction of $\mathcal{G}(\mathcal{T}_{t,n})$ from its point-process*].

Informally, think of forming this point-process by taking the heights of the branching points of the genealogical tree $\mathcal{G}(\mathcal{T}_{t,n})$ in the order they have as vertices in the tree. For convenience (in considering asymptotics with $t$ increasing) we keep track of the heights of the branching points in terms of their distances from level $t$. The vertical coordinate of each branching point is thus its distance below level $t$, while its horizontal coordinate is just its index. The point-process representation of $\mathcal{G}(\mathcal{T}_{t,n})$ from Figure 2(a) is shown in Figure 2(b). Formally, let $a_i, 1 \leq i \leq n-1$, be the times (distance to the root) of branch-points in the tree $\mathcal{G}(\mathcal{T}_{t,n})$, indexed in order induced from the depth-first search of the vertices in $\mathcal{T}_{t,n}$, let $t_i = t - a_i$ be their distance below level $t$ and let $\ell_i = i$.

DEFINITION 2. The *genealogical point-process* $\Pi_{t,n}$ is the random finite set

(1) $$\Pi_{t,n} = \{(\ell_i, t_i) : 1 \leq i \leq n-1, 0 < t_i < t\}.$$

For practical purposes it is most useful to exploit the bijection between a random tree and its contour process. We can obtain the point-process $\Pi_{t,n}$ equivalently from the contour process $\mathcal{C}_{\mathcal{T}_{t,n}}$ as follows. The $i$th individual extant at $t$ corresponds to the pair $(U_i, D_i)$: of the $i$th up-crossing time $U_i$ of level $t$ by the contour process $\mathcal{C}_{\mathcal{T}_{t,n}}$ and of the $i$th down-crossing time $D_i$ of level $t$. A precise definition of $U_i$ and $D_i$ will be given by (3) and (4) in the proof of Lemma 3. The branch-points $a_i, 1 \le i \le n-1$, of $\mathcal{G}(\mathcal{T}_{t,n})$ then correspond to the levels of lowest local minima of the excursions of $\mathcal{C}_{\mathcal{T}_{t,n}}$ below level $t$, in other words $a_i = \inf\{\mathcal{C}_{\mathcal{T}_{t,n}}(u) : D_i < u < U_{i+1}\}$.

We next use this observation together with the description of the law of $\mathcal{C}_{\mathcal{T}_{t,n}}$ to obtain the law of $\Pi_{t,n}$. We first recall the result of Neveu, Pitman and Le Gall, regarding the law of the contour process $\mathcal{C}_{\mathcal{T}}$ of an unconditioned random tree $\mathcal{T}$ (one can consult either [12] or [14] for its proof).

LEMMA 1. *In the contour process $\mathcal{C}_{\mathcal{T}}$ of a critical branching process $\mathcal{T}$ the sequence of rises and falls (up to the last fall) has the same distribution as a sequence of independent* Exponential($rate$ 1) *variables stopped one step before the sum of successive rises and falls becomes negative (the last fall is then set to equal this sum).*

The following corollary is an immediate consequence of Lemma 1 and the memoryless property of the exponential distribution.

COROLLARY 2. *For the contour process $\mathcal{C}_{\mathcal{T}}$ the process $X_{\mathcal{T}} = (\mathcal{C}_{\mathcal{T}}, \text{slope}[\mathcal{C}_{\mathcal{T}}])$ is a time-homogeneous strong Markov process on $\mathbb{R}^+ \times \{+1, -1\}$ stopped when it first reaches $(0, -1)$.*

The law of the genealogical point-process $\Pi_{t,n}$ can now easily be derived using some standard excursion theory of Markov processes. Note that the contour process of a whole class of binary branching processes can be shown to be a time-homogeneous Markov process as well (see [9]). In the following lemma we show that the distances of the $n-1$ branching points below level $t$ are independent and identically distributed, with the same law as that of the height of a random tree $\mathcal{T}$ conditioned on its height being less than $t$.

LEMMA 3. *For any fixed $t > 0$, the random set $\Pi_{t,n}$ is a simple point-process on $\{1, \ldots, n-1\} \times (0, t)$ with intensity measure*

(2) $$\nu_{t,n}(\{i\} \times d\tau) = \frac{d\tau}{(1+\tau)^2} \frac{1+t}{t}.$$

*In other words, $t_i$, $1 \le i \le n-1$, are i.i.d. variables on $(0, t)$ with the law* (2).

PROOF. In short the proof relies on the following. The contour process $\mathcal{C}_{\mathcal{T}}$ of an unconditioned tree $\mathcal{T}$ is, by the previous corollary, a Markov process considered until a certain stopping time. Hence, its excursions below some level $t$ are independent and identically distributed. Conditioning of the tree $\mathcal{T}_{t,n}$ translates simply in terms of its contour process, into conditioning this Markov process to have exactly $n-1$ excursions below $t$ until this stopping time. Further, for the law of these excursions it will follow, by the sign invariance of the law of $\mathcal{C}_\mathcal{T}$, that their law is the same as that of a copy of $\mathcal{C}_\mathcal{T}$ conditioned to have a height less than $t$.

Consider the Markov process $X_\mathcal{T} = (\mathcal{C}_\mathcal{T}, \text{slope}[\mathcal{C}_\mathcal{T}])$ until the first hitting time $U_{(0,-1)} = \inf\{u \geq 0 : X_\mathcal{T}(u) = (0,-1)\}$, and consider its excursions from the point $(t,+1)$ using the distribution of $\mathcal{C}_\mathcal{T}$ given by Lemma 1. For $i \geq 1$ let $U_i$ be the times of the up-crossings of level $t$ by $\mathcal{C}_\mathcal{T}$,

$$(3) \qquad U_0 = 0, \qquad U_i = \inf\{u > U_{i-1} : X_\mathcal{T}(u) = (t,+1)\}, \qquad i \geq 1.$$

Clearly $\mathbf{P}_{(t,+1)}[\inf\{u > 0 : X_\mathcal{T}(u) = (t,+1)\} > 0] = 1$; hence the set of all visits to $(t,+1)$ at times $\{U_i, i \geq 1\}$ is discrete. The excursions of $X_\mathcal{T}$ from level $t$ are, for $i \geq 1$,

$$e_i(u) = \begin{cases} X_\mathcal{T}(U_i + u), & \text{for } u \in [0, U_{i+1} - U_i), \\ (0,+1), & \text{else.} \end{cases}$$

The number of visits in an interval $[0,u]$ is

$$\ell(0) = 0, \qquad \ell(u) = \sup\{i > 0 : u > U_i\}, \qquad u > 0,$$

and the total number prior to $U_{(0,-1)}$ is $L = \sup\{i \geq 0 : U_{(0,-1)} > U_i\} = \ell(U_{(0,-1)})$. If $\mathbf{n}$ is the $\mathbf{P}_{(t,+1)}$-law of $e_i$, and if $\mathcal{E}^{<t}$ is the set of excursions from $(t,+1)$ that return to $(t,+1)$ without reaching $(0,-1)$, and $\mathcal{E}^{>t}$ the set of all others, then it is clear that (e.g., [19], **2**, Section VI.50) the following hold:

(a) $\mathbf{P}_{(t,+1)}[L \geq i] = [\mathbf{n}(\mathcal{E}^{<t})]^{i-1}$, $i \geq 1$, and $e_1, e_2, \ldots$ are independent;
(b) given that $L \geq i$, the law of $e_1, e_2, \ldots, e_{i-1}$ is $\mathbf{n}(\cdot \cap \mathcal{E}^{<t})/\mathbf{n}(\mathcal{E}^{<t})$;
(c) given that $L = i$, the law of $e_i$ is $\mathbf{n}(\cdot \cap \mathcal{E}^{>t})/\mathbf{n}(\mathcal{E}^{>t})$.

This makes $\{(\ell(U_i), e_i), 1 \leq i \leq L-1\}$ a simple point-process [note that $\ell(U_i) = i$, and $\ell(\infty) = L$], with the number of points having a Geometric($\mathbf{n}(\mathcal{E}^{>t})$) law, and with each $e_i$ having the law $\mathbf{n}(\cdot \cap \mathcal{E}^{<t})/\mathbf{n}(\mathcal{E}^{<t})$.

This observation is particularly convenient for analyzing the law of $\mathcal{C}_{\mathcal{T}_{t,n}}$. Since $\mathcal{C}_{\mathcal{T}_{t,n}}$ is just $\mathcal{C}_\mathcal{T}$ conditioned on $L = n$, the $n-1$ excursions of $\mathcal{C}_{\mathcal{T}_{t,n}}$ below $t$ are independent identically distributed with the law $\mathbf{n}(\cdot \cap \mathcal{E}^{<t})/\mathbf{n}(\mathcal{E}^{<t})$. We next derive the law of their depth $a_i$ measured as distance from level $t$ by $t_i = t - a_i$.

For each up-crossing time $U_i$ of level $t$, we have a down-crossing time

$$(4) \qquad D_i = \inf\{u > U_i : X_\mathcal{T}(u) = (t,-1)\}, \qquad i \geq 1.$$

For the values of $a_i$, $i \geq 1$, we are only interested in the part of the excursions from $(t,+1)$ below level $t$,

$$e_i^{<t} = e_i(D_i + u), \qquad u \in [0, U_{i+1} - D_i) \quad \text{and} \quad e_i^{<t}(u) = (0,+1) \qquad \text{else.}$$

We note that the shift and reflection invariance of the transition function of $\mathcal{C}_\mathcal{T}$, as well as its strong Markov property, applied to the law $\mathbf{n}$ for $e_i^{<t}$ imply that the law of $e_i^+ = t - e_i^{<t}$ is the same as the law of $X_\mathcal{T}$ conditioned to return to $(0, -1)$ before reaching $(t, +1)$. Consequently the law of $t - \inf(e_i^{<t}) = \sup(e_i^+)$ is the same law as that of $\sup(\mathcal{C}_\mathcal{T})$ conditioned to be less than $t$.

To explicitly express the law of $\sup(\mathcal{C}_\mathcal{T})$ we now recall classical results for the branching process $\mathcal{T}$ (e.g., [8], Section XVII.10.11), by which the law of the population size $N(t)$ of $\mathcal{T}$ at time $t$ is given by

$$(5) \quad \mathbf{P}[N(t) = 0] = \frac{t}{1+t}; \qquad \mathbf{P}[N(t) = k] = \frac{t^{k-1}}{(1+t)^{k+1}} \qquad \text{for } k \geq 1.$$

Hence

$$(6) \qquad \mathbf{P}[\sup(\mathcal{C}_\mathcal{T}) > t] = \mathbf{P}[N(t) > 0] = \frac{1}{1+t} \qquad \text{for } t \geq 0.$$

Now for $\mathcal{C}_{\mathcal{T}_{t,n}}$ and for each $1 \leq i \leq n-1$ we have that $a_i = \inf(e_i^{<t})$, and the $e_i^{<t}$ are independent with $e_i^{<t} \sim \mathbf{n}(\cdot \cap \boldsymbol{\mathcal{E}}^{<t})/\mathbf{n}(\boldsymbol{\mathcal{E}}^{<t})$; hence each $t_i = t - a_i$ has the law

$$(7) \quad \begin{aligned} \mathbf{P}[t_i \in d\tau] &= \mathbf{P}[\sup(\mathcal{C}_\mathcal{T}) \in d\tau | \sup(\mathcal{C}_\mathcal{T}) < t] \\ &= \frac{d\tau}{(1+\tau)^2} \frac{1+t}{t} \qquad \text{for } 0 \leq \tau \leq t. \end{aligned} \qquad \square$$

Asymptotics for $\Pi_{t,n}$ could now be established with a routine calculation. Instead of considering this result in isolation, it is far more natural to view it as part of the larger picture connecting critical branching processes and Brownian excursions. Let us recall the asymptotic results for critical Galton–Watson processes conditioned on a "large" total population size. A result of Aldous ([3], Theorem 23) says that its contour process (when appropriately rescaled) converges, as the total population size increases, to a Brownian excursion (doubled in height) conditioned to be of length 1. Note that, if $N_{\text{tot}}$ is the total population size of a critical Galton–Watson process, and $N(t_n)$ its population size at some given time $t_n$, then the events $\{N_{\text{tot}} = n\}$ and $\{N(t_n) = n | N(t_n) > 0\}$ are both events of "small" probabilities. The first has asymptotic chance $cn^{-3/2}$ as $n \to \infty$, and for $\{t_n\}_{n \geq 1}$ such that $t_n/n \to t$ as $n \to \infty$ the second has asymptotic chance $c(t)n^{-1}$ [3]. While the total population $N_{\text{tot}}$ size corresponds to the total length of the contour process, the population size $N(t_n)$ at a particular time $t_n$ corresponds to the occupation time of the contour process at level $t_n$. Hence, it is natural to expect that the contour process of a critical Galton–Watson process conditioned on a "large" population at time $t_n$ (when appropriately rescaled) converges, when $t_n/n \to t$ as $n \to \infty$, to a Brownian excursion conditioned to have local time 1 at level $t$.

We will show the following. Consider a Brownian excursion conditioned to have local time 1 at level $t$, as a "contour process" of an infinite tree (in the sense of the bijection between continuous functions and trees established in [3]). Consider defining a "genealogical" point-process from this Brownian excursion, using the depths of its excursions below level $t$, in the same manner as used in defining $\Pi_{t,n}$ from the contour process $\mathcal{C}_{\mathcal{T}_{t,n}}$, except that the excursions are now indexed by the amount of local time at level $t$ at their beginning. The state-space of such a point-process can be simply described, and we show that it has quite a simple law as well. It is then easy to show that this point-process is precisely the asymptotic process of appropriately rescaled processes $\Pi_{t_n,n}$ as $n \to \infty$.

We construct a point-process from a Brownian excursion conditioned to have local time 1 at level $t$, in the same manner in which $\Pi_{t,n}$ was constructed from the contour process $\mathcal{C}_{\mathcal{T}_{t,n}}$. Let $\mathcal{B}(u)$, $u \geq 0$, be a Brownian excursion. For a fixed $t > 0$, let $\ell_t(u)$, $u \geq 0$, be its local time at level $t$ up to time $u$ with the normalization of local time as one-half the occupation density relative to Lebesgue measure (the normalization choice is analogous to the upcrossings-only count for the contour process $\mathcal{C}_{\mathcal{T}}$). Let $i_t(\ell)$, $\ell \geq 0$, be the inverse process of $\ell_t$, in other words $i_t(\ell) = \inf\{u > 0 : \ell_t(u) > \ell\}$. Let $\mathcal{B}_{t,1}(u)$, $u \geq 0$, then be the excursion $\mathcal{B}$ conditioned to have total local time $\ell_t$ equal to 1, where $\ell_t = \ell_t(\infty)$ is the total local time at $t$. Consider excursions $e_\ell^{<t}$ of $\mathcal{B}_{t,1}$ below level $t$ indexed by the amount of local time $\ell$ at the time $i_t(\ell^-)$ of their beginning. For each such excursion let $a_\ell$ be its infimum, and let $t_\ell$ be the depth of the excursion measured from level $t$, $t_\ell = t - a_\ell$. Itô's excursion theory then insures that the process $\{(\ell, t_\ell) : i_t(\ell^-) \neq i_t(\ell)\}$ is well defined.

DEFINITION 3. The *continuum genealogical point-process* $\pi_{t,1}$ is the random countably infinite set

(8) $$\pi_{t,1} = \{(\ell, t_\ell) : i_t(\ell^-) \neq i_t(\ell)\}.$$

REMARK 4. The name of the process will be justified by establishing it as the limit of genealogical point-processes.

For the state-space of the continuum genealogical process we introduce the notion of a nice point-process ([3], Section 2.8). A *nice point-process on* $[0,1] \times (0, \infty)$ is a countably infinite set of points such that the following hold:

1. for any $\delta > 0$, $[0,1] \times [\delta, \infty)$ contains only finitely many points;
2. for any $0 \leq x < y \leq 1$, $\delta > 0$, $[x,y] \times (0, \delta)$ contains at least one point.

We next show that the state-space for $\pi_{t,1}$ is the set of nice point-processes, and we establish the law of this process using standard results of Levy, Itô and Williams on excursion theory.

LEMMA 4. *The random set $\pi_{t,1}$ is a Poisson point-process on $[0,1] \times (0,t)$ with intensity measure*

$$\nu(d\ell \times d\tau) = d\ell \frac{d\tau}{\tau^2}. \tag{9}$$

*In particular, the random set $\pi_{t,1}$ is a.s. a nice point-process.*

PROOF. The crux of the proof lies in the following observations. An unconditioned Brownian excursion $\mathcal{B}$ observed from the first time it reaches level $t$ is just $t$—a standard Brownian motion observed until the first time it reaches $t$. The excursions of $\mathcal{B}$ below level $t$ are thus the positive excursions of the Brownian motion. By a standard result, the process of excursions of Brownian motion from 0, indexed by the amount of local time at 0 at the time of their beginning, is a Poisson point-process with intensity measure $d\ell \times \mathbf{n}$, where $\mathbf{n}$ is Itô's excursion measure. One can show that the condition on $\mathcal{B}$ to have local time 1 at level $t$ is equivalent to the condition that the shifted Brownian motion has all its excursions until local time 1 of height lower than $t$ and has one excursion at local time 1 higher than $t$. This then, by the independence properties of Poisson processes, allows for a simple description of the point-process of the depths of excursions below $t$ of $\mathcal{B}_{t,1}$ as a Poisson process itself, except restricted to the set $[0,1] \times (0,t)$.

Consider the path of an (unconditioned) Brownian excursion $\mathcal{B}$ after the first hitting time of $t$, $U_t = \inf\{u \geq 0 : \mathcal{B}(u) = t\}$, shifted and reflected about the $u$-axis

$$\beta(u) = t - \mathcal{B}(U_t + u) \quad \text{for } u \geq 0. \tag{10}$$

Let $\ell_0^\beta(u)$, $u \geq 0$, be the local time of $\beta$ at level 0 up to time $u$, and let $i_0^\beta(\ell)$, $\ell \geq 0$, be the inverse process of this local time, in other words $i_0^\beta(\ell) = \inf\{u > 0 : \ell_0^\beta(u) > \ell\}$. Then the process $\beta(u)$, $u \geq 0$, is a standard Brownian motion stopped at the first hitting time of $t$, $U_t^\beta = \inf\{u \geq 0 : \beta(u) = t\}$.

Next, the excursions of $\beta$ from 0 are (with a change of sign) precisely the excursions of $\mathcal{B}$ from $t$, and the local time process $\ell_0^\beta$ of $\beta$ is equivalent to the local time process $\ell_t$ of $\mathcal{B}$. We are only interested in the excursions of $\mathcal{B}$ below $t$, which are the positive excursions of $\beta$ indexed by $\ell$ such that $i_0^\beta(\ell^-) \neq i_0^\beta(\ell)$ and $\sup\{\beta(u) : i_0^\beta(\ell^-) \leq u \leq i_0^\beta(\ell)\} > 0$,

$$e_\ell^+ = \beta(i_0^\beta(\ell^-) + u), \quad u \in [0, i_0^\beta(\ell) - i_0^\beta(\ell^-)) \quad \text{and} \quad e_\ell^+(u) = 0 \quad \text{else.}$$

Note that we thus have that the infimum of an excursion of $\mathcal{B}$ below $t$ to be simply $\inf(e_\ell^{<t}) = t - \sup(e_\ell^+)$.

Standard results of Itô's excursion theory (e.g., [19], **2**, Section VI.47) imply that for a standard Brownian motion $\beta$ the random set $\{(\ell, \sup(e_\ell^+)) : i_0^\beta(\ell^-) \neq$ $i_0^\beta(\ell)\}$ is a Poisson point-process on $\mathbb{R}^+ \times \mathbb{R}^+$ with intensity measure $d\ell\, d\tau/\tau^2$ (recall our choice for the normalization of local time).

Now let $\mathrm{L} = \inf\{\ell \geq 0 : \sup(e_\ell^+) \geq t\}$. Then stopping $\ell_0^\beta$ at the hitting time L is equivalent to stopping $\beta$ at its hitting time $U_t^\beta$. Let $\pi_t$ be a random set defined from the unconditioned Brownian excursion $\mathcal{B}$, in the same manner in which we defined $\pi_{t,1}$ from a conditioned Brownian excursion $\mathcal{B}_{t,1}$. Then, using the relationship (10) of $\mathcal{B}$ and $\beta$, we observe that $\pi_t$ is equivalent to a restriction of $\{(\ell, \sup(e_\ell^+)) : i_0^\beta(\ell^-) \neq i_0^\beta(\ell)\}$ on the random set $[0, \mathrm{L}] \times (0, t)$. The Poisson point-process description of $\{(\ell, \sup(e_\ell^+)) : i_0^\beta(\ell^-) \neq i_0^\beta(\ell)\}$ now implies that $\pi_t$ is a Poisson point-process on $\mathbb{R}^+ \times \mathbb{R}^+$ with intensity measure $d\ell\, d\tau/\tau^2$ restricted to the random set $[0, \mathrm{L}] \times (0, t)$.

Next, note that the condition $\{\ell_t = 1\}$ for $\mathcal{B}$ is equivalent to the condition $\{\ell_0^\beta(U_t^\beta) = 1\}$ for $\beta$, which is further equivalent to the condition $\{\mathrm{L} = 1\}$ for $\pi_t$. We have thus established that $\pi_{t,1} \stackrel{d}{=} \pi_t | \{\mathrm{L} = 1\}$.

Further, the condition $\{\mathrm{L} = 1\}$ on $\pi_t$ is equivalent to the condition that $\pi_t$ has no points in $[0, 1) \times [t, \infty)$ and has a point in $\{1\} \times [t, \infty)$. However, since $\pi_t$ is Poisson, independence of Poisson random measures on disjoint sets implies that conditioning $\pi_t$ on $\{\mathrm{L} = 1\}$ will not alter its law on the set $[0, 1] \times (0, t)$. However, since $\pi_{t,1}$ is supported precisely on $[0, 1] \times (0, t)$, the above results together imply that $\pi_{t,1}$ is a Poisson point-process on $[0, 1] \times (0, t)$ with intensity measure $d\ell\, d\tau/\tau^2$.

It is now easy to see from the intensity measure of $\pi_{t,1}$ that its realizations are a.s. nice point-processes, namely:

(a) for any $\delta > 0$, $\iint_{[0,1]\times[\delta,\infty)} d\ell\, d\tau/\tau^2 = 1/\delta < \infty$;
(b) for any $0 \leq x < y \leq 1$, and $\delta > 0$, $\iint_{[x,y]\times(0,\delta)} d\ell\, d\tau/\tau^2 = (y-x) \cdot \infty$.

Also, since $\pi_{t,1}$ is Poisson, finiteness of its intensity measure on $[0, 1] \times [\delta, \infty)$ implies that it has a.s. only finitely many points in the set $[0, 1] \times [\delta, \infty)$, while infiniteness of its intensity measure on $[x, y] \times (0, \delta)$ implies that it has a.s. at least one point on the set $[x, y] \times (0, \delta)$. □

Having thus obtained the description of the continuum genealogical point-process induced by a conditioned Brownian excursion, it is now a simple task to confirm that it indeed arises as the limit of genealogical processes. The right rescaling for $\mathcal{T}_{t,n}$ is to speed up the time by $n$ and to assign mass $n^{-1}$ to each extant individual, which implies the appropriate rescaling of each coordinate of $\Pi_{t,n}$ by $n^{-1}$. We hence define the rescaled genealogical point-process as

$$(11) \qquad n^{-1}\Pi_{t,n} = \{(n^{-1}\ell_i, n^{-1}t_i) : (\ell_i, t_i) \in \Pi_{t,n}\}$$

and establish its asymptotic behavior as $n \to \infty$.

THEOREM 5. *For any $\{t_n > 0\}_{n \geq 1}$ such that $t_n/n \xrightarrow[n \to \infty]{} t$ we have*

$$n^{-1} \Pi_{t_n,n} \xRightarrow[n \to \infty]{d} \pi_{t,1}.$$

REMARK 5. The notation $\xRightarrow{d}$ is used to mean weak convergence of processes.

PROOF OF THEOREM 5. The proof of the theorem is a just consequence of the fact that weak convergence of Poisson point-processes follows from the weak convergence of their intensity measures.

By Lemma 3 and the rescaling (11) we have that $n^{-1} \Pi_{t_n,n}$ is a simple point-process on $\{1/n, \ldots, 1-1/n\} \times (0, t_n/n)$ with intensity measure

$$\text{(12)} \qquad \frac{1}{n} \sum_{i=1}^{n-1} \delta_{\{i/n\}}(\ell) \frac{n \, d\tau}{(1+n\tau)^2} \frac{1+t_n}{t_n}.$$

If $\{t_n\}_{n \geq 1}$ is such that $t_n/n \to t$ as $n \to \infty$, then it is clear that the support set of the process $n^{-1} \Pi_{t_n,n}$ converges to $[0,1] \times (0,t)$, the support set of the process $\pi_{t,1}$. It is also clear that the intensity measure (12) converges to $d\ell \, d\tau / \tau^2$, which, by Lemma 4, is the intensity measure of $\pi_{t,1}$. For simple point-processes this is sufficient (e.g., [5], Section 12.3) to insure weak convergence of the processes $n^{-1} \Pi_{t_n,n}$ to a Poisson point-process on $[0,1] \times (0,t)$ with intensity measure $d\ell \, d\tau / \tau^2$. By Lemma 4, we thus have that $n^{-1} \Pi_{t_n,n} \xRightarrow[n \to \infty]{d} \pi_{t,1}$. □

**3. Genealogy of sampled extinct individuals.** We now want to extend the analysis of the ancestry of extant individuals to include some proportion of the extinct individuals as well. Suppose that each individual in the past has independently had a chance $p$ of appearing in the historical record. We indicate such sampling of extinct individuals by putting a star mark on the leaf of $\mathcal{T}_{t,n}$ corresponding to the recorded individual. An example of a realization of such $p$-sampling is shown in Figure 3(a), and the induced sampling in the contour process is shown in Figure 3(b).

The goal is to combine the information on the sampled extinct individuals with our analysis of the ancestry of the extant ones. To do so we extend our earlier notions of the genealogy of the extant individuals and of the genealogical point-process.

Let the *p-sampled history* of extant individuals at time $t$ be defined as the smallest subtree of the family tree which contains all the edges representing the ancestry both of the extant individuals and of all of the $p$-sampled extinct individuals. We denote the $p$-sampled history of extant individuals at $t$ in $\mathcal{T}_{t,n}$ by $\mathcal{G}_p(\mathcal{T}_{t,n})$. Note that by definition $\mathcal{G}_p(\mathcal{T}_{t,n})$ contains the genealogy $\mathcal{G}(\mathcal{T}_{t,n})$ (which would correspond to a 0-sampled history). It is in fact convenient to think of $\mathcal{G}_p(\mathcal{T}_{t,n})$ as consisting of the "main genealogical tree" $\mathcal{G}(\mathcal{T}_{t,n})$ and a collection of "$p$-sampled subtrees" attached to this main tree linking with additional branches the ancestry of $p$-sampled extinct individuals. Figure 4(a) shows the $p$-genealogical subtree of the tree from Figure 3(a). We next extend the notion of the genealogical point-process to represent this enriched $p$-sampled genealogy. We construct a point-process representation of $\mathcal{G}_p(\mathcal{T}_{t,n})$ so that it contains $\Pi_{t,n}$ as its "main points."

Informally, think of extending the point-process $\Pi_{t,n}$ [representing $\mathcal{G}(\mathcal{T}_{t,n})$], by adding sets representing the $p$-sampled subtrees as follows. At each branch-point of the main tree there is a set of $p$-sampled subtrees attached to the edges of the main tree on the left of this branching point, and a set of $p$-sampled subtrees attached on the right of this branching point [see Figure 4(a)]. We associate with each branch-point at height $a_i$ a left set $\mathcal{L}_i$ and a right set $\mathcal{R}_i$, which shall represent these sets of subtrees. Each such $\mathcal{L}_i$ and $\mathcal{R}_i$ needs to contain the following information: the heights $a_{i,L}(j)$ and $a_{i,R}(j)$ at which the $p$-sampled subtrees get attached to the edges of the main tree [as before we shall keep track of these heights as distances from level $t$ in terms of $t_{i,L}(j) = t - a_{i,L}(j)$ and $t_{i,R}(j) = t - a_{i,R}(j)$]; and the shape of the subtrees $\Upsilon_{i,L}(j)$ and $\Upsilon_{i,R}(j)$ themselves (the indexing $j \geq 0$ on the subtrees is induced by a depth-first search forward to the branch-point at $a_i$ for the left sets and a depth-first search backward to the branch-point at $a_i$ for the right sets). The point-process representing the $p$-sampled genealogical tree from Figure 4(a) is shown in Figure 4(b). To describe the law of the $p$-subtrees it will also be convenient to keep track of the height $h_{i,L}(j)$ and $h_{i,R}(j)$ of the subtrees $\Upsilon_{i,L}(j)$ and $\Upsilon_{i,R}(j)$.

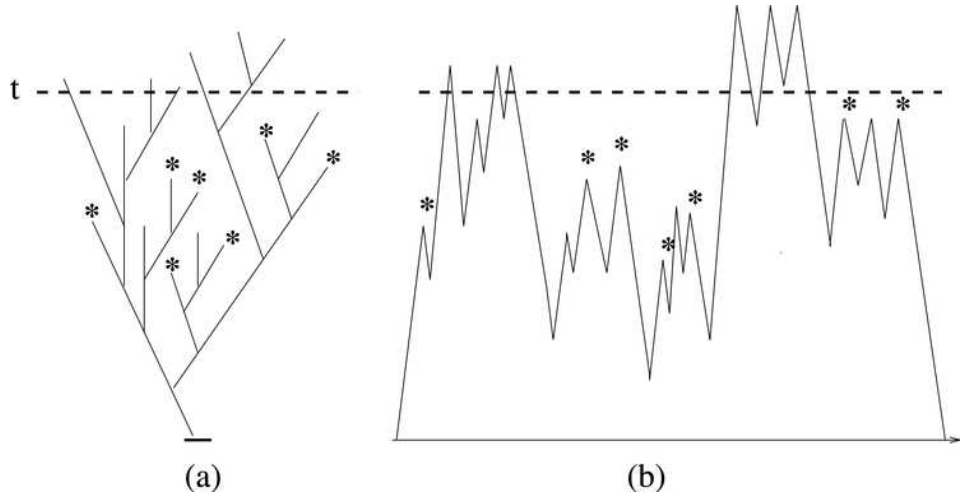

FIG. 3. (a) *The tree $\mathcal{T}_{t,n}$ with $p$-sampling on its individuals (the sampled individuals are represented by $*$'s).* (b) *The contour process of this tree with the sampling on the corresponding local maxima.*

Formally, we define the point-process of $\mathcal{G}_p(\mathcal{T}_{t,n})$ from the contour process $\mathcal{C}_{\mathcal{T}_{t,n}}$. The $p$-sampling on the tree is represented by the sampling of the local maxima of $\mathcal{C}_{\mathcal{T}_{t,n}}$. From the definition of $\Pi_{t,n}$, we have the heights of the branch-points of $\mathcal{G}(\mathcal{T}_{t,n})$ to be $a_i = \inf\{\mathcal{C}_{\mathcal{T}_{t,n}}(u) : D_i < u < U_{i+1}\}$, occurring in the contour process $\mathcal{C}_{\mathcal{T}_{t,n}}$ at times $B_i = \arg\min\{\mathcal{C}_{\mathcal{T}_{t,n}}(u) : u \in (D_i, U_{i+1})\}$. The set $\mathcal{L}_i$, representing the set of $p$-subtrees attaching to the edges of $\mathcal{G}(\mathcal{T}_{t,n})$ on the left of the branch-point $a_i$, is defined from the part of the excursion of $\mathcal{C}_{\mathcal{T}_{t,n}}$ below $t$ before time $B_i$. In other words if, for $X_{\mathcal{T}_{t,n}} = (\mathcal{C}_{\mathcal{T}_{t,n}}, \text{slope}[\mathcal{C}_{\mathcal{T}_{t,n}}])$, we define

$$e_{i,L}^{<t}(u) = X_{\mathcal{T}_{t,n}}(D_i + u), \qquad u \in [0, B_i - D_i),$$

then $\mathcal{L}_i$ is completely defined by $e_{i,L}^{<t}$. Analogously $\mathcal{R}_i$ is defined from the part of the excursion of $\mathcal{C}_{\mathcal{T}_{t,n}}$ below $t$ after time $B_i$; in other words if we define

$$e_{i,R}^{<t}(u) = X_{\mathcal{T}_{t,n}}(U_{i+1} - u), \qquad u \in [0, U_{i+1} - B_i),$$

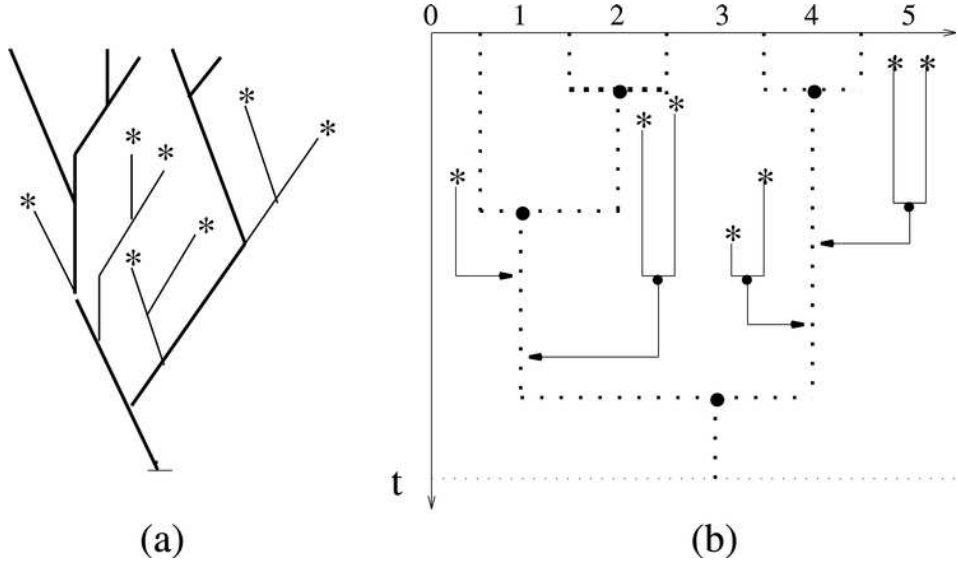

FIG. 4. (a) *The $p$-sampled tree $\mathcal{G}_p(\mathcal{T}_{t,n})$; the "main tree" (in bold) has the "$p$-sampled subtrees" attached to it.* (b) *The point-process representation $\Xi_{t,n}^p$ of $\mathcal{G}_p(\mathcal{T}_{t,n})$; each of the "main points" (large dots) has an associated left set and a right set representing the $p$-sampled subtrees attaching to the left and right of that branch-point.*

then it is completely defined by $e_{i,R}^{<t}$ (the subscripts $L$ and $R$ reflect whether the entities are involved in defining $\mathcal{L}_i$ or $\mathcal{R}_i$). Note that the $e_{i,L}^{<t}$ runs forward up to time $B_i$, whereas $e_{i,R}^{<t}$ runs backward. On the extreme ends, we have the set of $p$-subtrees on the left of the first branching point defined by the part of $\mathcal{C}_{\mathcal{T}_{t,n}}$ prior to the first up-crossing time $U_1$, $e_{0,R}^{<t}(u) = X_{\mathcal{T}_{t,n}}(U_1 - u)$, $u \in [0, U_1)$. Analogously, the set of $p$-subtrees on the right of the last branching point is defined by the part of $\mathcal{C}_{\mathcal{T}_{t,n}}$ after the last down-crossing time $D_n$, $e_{n,L}^{<t}(u) = X_{\mathcal{T}_{t,n}}(D_n + u)$, $u \in [0, U_{(0,-1)} - D_n)$, where $U_{(0,1)} = \inf\{u \geq 0 : X_{\mathcal{T}_{t,n}} = (0,-1)\}$.

To define the sets $\mathcal{L}_i$ and $\mathcal{R}_i$ we also need to define the processes

$$\varsigma_{i,L}(u) = \inf_{0 \leq v \leq u} e_{i,L}^{<t}(v), \qquad u \in [0, B_i - D_i),$$

$$\varsigma_{i,R}(u) = \inf_{0 \leq v \leq u} e_{i,R}^{<t}(v), \qquad u \in [0, U_{i+1} - B_i).$$

The bijection between the tree $\mathcal{T}_{t,n}$ and its contour process $\mathcal{C}_{\mathcal{T}_{t,n}}$ implies that the heights at which the $p$-subtrees are attached to the edges of the main tree are precisely the levels of constancy of the processes $\varsigma_{i,L}$ and $\varsigma_{i,R}$. Furthermore, the $p$-subtrees themselves have as their contour processes the excursions of $e_{i,L}^{<t} - \varsigma_{i,L}$ and $e_{i,R}^{<t} - \varsigma_{i,R}$ above these levels of constancy (see [16] for a detailed description). Figure 5 shows $e_{i,L}^{<t}$ together with its infimum process $\varsigma_{i,L}^{<t}$.

We define $a_{i,L}(j)$, $j \geq 0$, to be the successive levels of constancy of $\varsigma_{i,L}$, and let $t_{i,L}(j) = t - a_{i,L}(j)$ be their distance from level $t$. For each level of constancy $a_{i,L}(j)$, let $e_{i,L}^{<t}(j)$ be the excursion of $e_{i,L}^{<t} - \varsigma_{i,L}$ that lies above the level $a_{i,L}(j)$. Let $h_{i,L}(j)$ be the height of this excursion, $h_{i,L}(j) = \sup(e_{i,L}^{<t}(j))$, and let $\Upsilon_{i,L}(j)$ be the tree whose contour process is the excursion $e_{i,L}^{<t}(j)$. Figure 5 shows an excursion $e_{i,L}^{<t}(j)$ with the $p$-subtree $\Upsilon_{i,L}(j)$ it defines. Note that all the star marks due to $p$-sampling are contained in the excursions $e_{i,L}^{<t}(j)$, hence are contained in the subtrees $\Upsilon_{i,L}(j)$. An analogous definition leads to $a_{i,R}(j)$, $j \geq 0$, $h_{i,R}(j)$ and $\Upsilon_{i,L}(j)$ from $e_{i,R}^{<t}(j)$ and $\varsigma_{i,R}(j)$. With each point $(\ell_i, t_i)$ of $\Pi_{t,n}$ we now associate the sets

(13)   $\mathcal{L}_i = \{(t_{i,L}(j), \Upsilon_{i,L}(j))\}_{j \geq 0}$  and  $\mathcal{R}_i = \{(t_{i,R}(j), \Upsilon_{i,R}(j))\}_{j \geq 0}$.

In addition, for extreme ends we define one set $\mathcal{R}_0$ from $e_{0,R}^{<t}$, and we define a set $\mathcal{L}_n$ from $e_{n,L}^{<t}$. For ease of future notation we set $\mathcal{L}_0 = \varnothing$, $\mathcal{R}_n = \varnothing$, $(\ell_0, t) = (1, t)$ and $(\ell_n, t_n) = (n, t)$.

DEFINITION 6.   The *$p$-sampled historical point-process* $\Xi_{t,n}^p$ is the random set

(14)         $\Xi_{t,n}^p = \{(\ell_i, t_i, \mathcal{L}_i, \mathcal{R}_i) : (\ell_i, t_i) \in \Pi_{t,n}, 0 \leq i \leq n\}.$

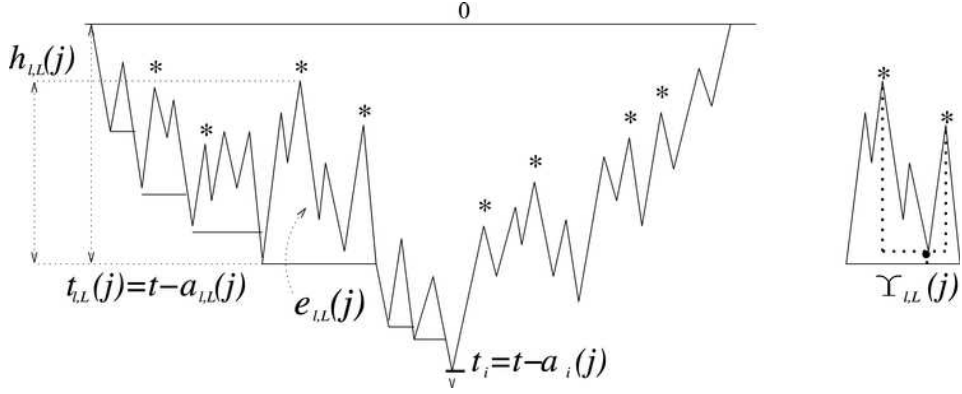

FIG. 5. *The left half $e_{i,L}^{<t}$ of an excursion of $\mathcal{C}_{\mathcal{T}_{t,n}}$ below $t$, with its infimum process $\varsigma_{i,L}$ whose levels of constancy are $\{a_{i,L}(j)\}_j$, above which lie the $p$-marked subtrees $\{\Upsilon_{i,L}(j)\}_j$ of heights $\{h_{i,L}(j)\}_j$.*

REMARK 7. We have in fact implicitly defined a point-process representation $\Xi_{t,n}$ of a complete historical point-process (which would correspond to 1-sampling). The difference between $\Xi_{t,n}$ and $\Xi_{t,n}^p$ is only in the $*$'s on the leaves in the latter. It will, however, be clear that for nice asymptotic behavior we need to consider $\Xi_{t,n}^p$ with $p < 1$; in other words we can only keep track of a proportion of the extinct individuals.

We can now derive the law of the point-process $\Xi_{t,n}^p$. For this we shall also need the law of the $p$-subtrees appearing in the sets $\mathcal{L}_i$ and $\mathcal{R}_i$. Let $\mathbf{T}$ denote the space of finite rooted binary trees with edge-lengths, and let $\Lambda$ denote the law on $\mathbf{T}$ of the tree $\mathcal{T}$. Then, let $\Lambda^p$ denote the law on $\mathbf{T}$ induced by the $p$-sampling on the tree $\mathcal{T}$. Further, for any $h > 0$, let $\Lambda_h^p$ denote the law induced by restricting $\Lambda^p$ to the trees $\mathcal{T}$ of height $h$.

To describe the law of $\Xi_{t,n}^p$ we use a more careful and detailed analysis of the structure of the contour process $\mathcal{C}_{\mathcal{T}_{t,n}}$. First we use the result of Lemma 3, which gives us the law of the main points of $\Xi_{t,n}^p$. Then conditional on the location of the main points, we give the law of the sets $\mathcal{L}_i$ and $\mathcal{R}_i$ of $p$-subtrees. We show that the sets $\mathcal{L}_i$ and $\mathcal{R}_i$ are independent Poisson point-processes. The intensity measure of each such set is given by the following. First, choose $t_{i,L}(j)$, the distances below $t$ at which the $p$-sampled subtrees are getting attached, uniformly over $t_i$, the total distance below $t$ to the $i$th branch-point. Next, choose $h_{i,L}(j)$, the height for each $p$-subtree, according to the same law as that of the height of a tree $\mathcal{T}$ whose height is known to be less than $t_{i,L}(j)$. Finally, choose the law of $\Upsilon_{i,L}(j)$, the attaching subtree, according to the law $\Lambda_h^p$ described above.

LEMMA 6. *For any fixed $0 < p < 1$, the law of the random set $\Xi_{t,n}^p$ is given by the following:*

(a) $\{(\ell_i, t_i) : 1 \leq i \leq n-1\}$ *is the simple point-process* $\Pi_{t,n}$ *of Lemma* 3*;*

(b) *given* $\{(\ell_i, t_i), 1 \leq i \leq n-1\}$*, the sets* $\mathcal{L}_i$ *and* $\mathcal{R}_i$ *are independent; and for each* $0 \leq i \leq n$ *the random sets* $\mathcal{L}_i$ *and* $\mathcal{R}_i$ *are Poisson point-processes on* $\mathbb{R}^+ \times \mathbf{T}$ *with intensity measure*

$$(15) \qquad \mathbb{1}_{\{0 < t < t_i\}} \, dt \, \mathbb{1}_{\{0 < h < t\}} \frac{dh}{(1+h)^2} \frac{1+t}{t} \Lambda_h^p.$$

PROOF. The proof relies on exploiting the alternating Exponential(rate 1) step structure of the contour process $\mathcal{C}_{\mathcal{T}_{t,n}}$. From Lemma 3 we have that the excursions of $\mathcal{C}_{\mathcal{T}_{t,n}}$ below $t$ are independent and that their law is the same as that of $\mathcal{C}_{\mathcal{T}}$ conditioned on having height less than $t$. We further show that, when decomposed into the part before its lowest point and a part after it, the two parts of these excursions are conditionally independent given the excursion's depth. In fact, if the former is run forward to the lowest point, and the latter backward to the lowest point, these two parts have the same conditioned law. We obtain a simple description of the law of the levels of constancy of the infimum process for these parts, and the excursions above these levels of constancy are shown to be copies of $\mathcal{C}_{\mathcal{T}}$ restricted in their height.

Independence of the pairs of sets $\mathcal{L}_i, \mathcal{R}_i$ over the index $i$ follows from the independence of the excursions $e_i^{<t}$ of $X_{\mathcal{T}_{t,n}}$ below level $t$ as shown in Lemma 3. The strong Markov property of $X_{\mathcal{T}}$ also gives the independence of $\mathcal{R}_0$ and $\mathcal{L}_n$ from all the pairs $\mathcal{L}_i, \mathcal{R}_i$. The proof of Lemma 3 also shows that the law of the excursions $t - e_i^{<t}$ is the same as that of $X_{\mathcal{T}}$ conditioned on $\sup(\mathcal{C}_{\mathcal{T}}) < t$. The left half $t - e_{i,L}^{<t}$ of such an excursion is defined as the part of $t - e_i^{<t}$ until it reaches its maximum, and the right half $t - e_{i,R}^{<t}$ as the part after this maximum, run backward in time (the $u$-coordinate). To derive the conditional law of $\mathcal{L}_i, \mathcal{R}_i$ given $t_i$ we thus need to analyze the conditional law of the two parts of $X_{\mathcal{T}}$ on either side of its maximum, given the maximum's value.

Let us first consider the process $X_{\mathcal{T}}$ continued past its first hitting time of $(0, -1)$. Let $S_{\mathcal{T}}$ be its maximum process

$$S_{\mathcal{T}}(u) = \sup_{0 \leq v \leq u} \mathcal{C}_{\mathcal{T}}(v), \qquad u \geq 0,$$

and let us consider the process $(S_{\mathcal{T}}, S_{\mathcal{T}} - \mathcal{C}_{\mathcal{T}})$, which clearly completely describes $X_{\mathcal{T}}$. The process $(S_{\mathcal{T}}, S_{\mathcal{T}} - \mathcal{C}_{\mathcal{T}})$ consists of an alternating sequence of the following: rises of slope 1 for $S_{\mathcal{T}}$ paired with the intervals at 0 for $S_{\mathcal{T}} - \mathcal{C}_{\mathcal{T}}$, and levels of constancy for $S_{\mathcal{T}}$ paired with the excursions from 0 for $S_{\mathcal{T}} - \mathcal{C}_{\mathcal{T}}$. Figure 6(a) shows a decomposition of $X_{\mathcal{T}}$ into $S_{\mathcal{T}}$ and $S_{\mathcal{T}} - \mathcal{C}_{\mathcal{T}}$. Because the alternating steps of $\mathcal{C}_{\mathcal{T}}$ are independent Exponential(rate 1) variables, it is not difficult to see that the alternating steps of $(S_{\mathcal{T}}, S_{\mathcal{T}} - \mathcal{C}_{\mathcal{T}})$ are independent, and that the rises of $S_\mathcal{T}$ all have the same Exponential(rate 1) distribution while the excursions of $S_\mathcal{T}-\mathcal{C}_\mathcal{T}$ all have the same distribution as $X_\mathcal{T}$ [stopped when it first hits $(0,-1)$]. Namely, the first rise of $S_\mathcal{T}$ is just the first rise of $\mathcal{C}_\mathcal{T}$; the $X_\mathcal{T}$ law and independence of a subsequent excursion of $S_\mathcal{T} - \mathcal{C}_\mathcal{T}$ is immediate from the law of $\mathcal{C}_\mathcal{T}$; and, finally, the Exponential(rate 1) law and independence of a subsequent rise of $S_\mathcal{T}$ is just a consequence of the memoryless property of rises of $\mathcal{C}_\mathcal{T}$.

Consider now the point-process $\{(s,\varepsilon_s)\}$ of levels of constancy of $S_\mathcal{T}$ paired with the excursions of $S_\mathcal{T}-\mathcal{C}_\mathcal{T}$ below them. The above analysis shows that $\{s(u),\ u\geq 0\}$ is a Poisson(rate 1) process, and the excursions $\varepsilon_s$ all have the same law as $\mathcal{C}_\mathcal{T}$. We shall denote the law of $\mathcal{C}_\mathcal{T}$ by $\mathbf{n}$ (as in the proof of Lemma 3). Then $\{(s,\varepsilon_s)\}$ forms a Poisson point-process with intensity measure $ds\,\mathbf{n}$. Note that it was shown in (6) that the height of $\mathcal{C}_\mathcal{T}$, and hence the height of these excursions, is given by $\mathbf{n}(\sup(\cdot)>h)=1/(1+h)$, for $h\geq 0$.

We now consider how the maximum of $X_\mathcal{T}$ in the time interval before it first hits $(0,-1)$ appears within this point-process $\{(s,\varepsilon_s)\}$. We denote the first hitting time of $(0,-1)$ by $U_{(0,-1)}=\inf\{t:X_\mathcal{T}=(0,-1)\}$, and the maximum of interest by $M=\sup\{X_\mathcal{T}(v):v\in(0,U_{(0,-1)})\}$. Then, one can easily note that $s(U_{(0,-1)})=M$ and that $\forall s\in(0,M)$, $h(\varepsilon_s)<s$ whereas for $s=M$, $h(\varepsilon_M)\geq M$. Figure 6(b) depicts a realization of $\{(s,\varepsilon_s)\}$. In other words, the point $(M,\varepsilon_M)$ is the first point (in terms of the $s$-coordinate) of the process $\{(s,\varepsilon_s)\}$ which falls outside the set $\{(s,h_s):s\geq 0,\ h_s\geq 0,\ h_s<s\}$. Independence of Poisson random measures on disjoint sets then implies that, given the value of $M=\sup\{X_\mathcal{T}(v):v\in(0,U_{(0,-1)})\}$, the conditional law of the process $\{(s,\varepsilon_s):s<M\}$ is independent of the point $(M,\varepsilon_M)$ and

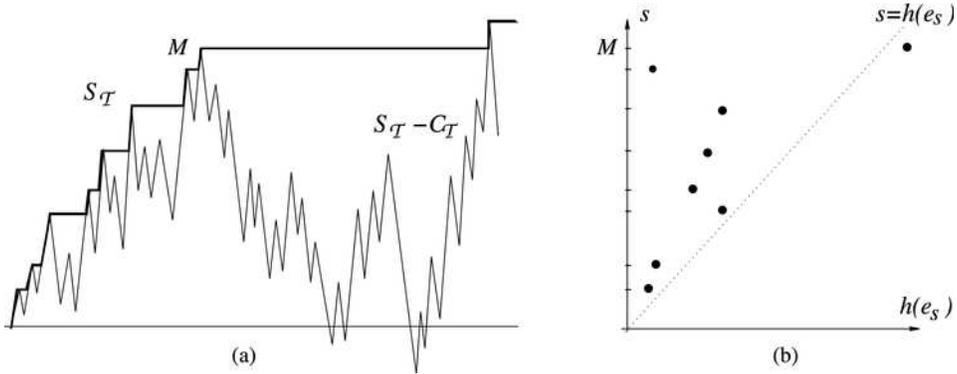

Fig. 6. (a) *The process $X_\mathcal{T}$ continued past its first hitting of $(0,-1)$, its maximum process $S_\mathcal{T}$ and the excursion process $S_\mathcal{T} - \mathcal{C}_\mathcal{T}$ below the levels of constancy of $S_\mathcal{T}$; $M$ is the maximal value of $X_\mathcal{T}$ before it first hits $(0,-1)$.* (b) *The point-process $\{(s,h(\varepsilon_s))\}$ of the values of constancy of $S_\mathcal{T}$, paired with the heights of excursions of $S_\mathcal{T}-\mathcal{C}_\mathcal{T}$ below them.*

has the law of a Poisson point-process with intensity measure

$$\mathbb{1}_{0<s<M}\,ds\,\mathbb{1}_{0<h<s}\frac{dh}{(1+h)^2}\mathbf{n}(\,\cdot\,|\sup(\cdot)=h).$$

It is important now to note that the process $\{(s,\varepsilon_s)\colon s<M\}$ completely describes the part of $X_\mathcal{T}$ before it reaches the maximum $M$, while the point $(M,\varepsilon_M)$ completely describes the part of $X_\mathcal{T}$ after it reaches $M$ and before it first hits $(0,-1)$.

We can now tend to the quantities of interest within the excursions $t-e_i^{<t}$. We have that the conditional law of $t-e_{i,L}^{<t}$, $t-e_{i,R}^{<t}$ given $t_i$ is the same as the conditional law of the parts of $X_\mathcal{T}$ before and after its maximum $M$ given that $M=t_i$. Within $t-e_{i,L}^{<t}$ the levels of constancy $t_{i,L}(j)$ and the associated excursions $t-e_{i,L}^{<t}(j)$ above these levels precisely correspond to the levels of constancy $s$ and its associated excursions $\varepsilon_s$ within the process $\{(s,\varepsilon_s)\colon s<t_i\}$. Moreover, $t-e_{i,R}^{<t}$ precisely corresponds to the part of the excursion $\varepsilon_M$ before it hits $(0,-1)$, reversed in time (cf. Figures 5 and 6).

From our analysis above it thus follows that $\{(t_{i,L}(j),e_{i,L}^{<t}(j))\}$, that is, $t-e_{i,L}^{<t}$ and $t-e_{i,R}^{<t}$ are conditionally independent given $t_i$, and that $\{(t_{i,L}(j),e_{i,L}^{<t}(j))\}$ has the law of a Poisson point-process with intensity measure

$$\mathbb{1}_{0<t<t_i}\,dt\,\mathbb{1}_{0<h<t}\frac{dh}{(1+h)^2}\mathbf{n}(\,\cdot\,|\sup(\cdot)=h).$$

Now the strong Markov property implies that the $p$-sampling on the local maxima of the whole contour process $\mathcal{C}_{\mathcal{T}_{t,n}}$ is for each $e_{i,L(j)}^{<t}$ again a Bernoulli $p$-sampling on its local maxima. Thus the conditional law of the $p$-sampled tree $\Upsilon_{i,L}(j)$ defined from $e_{i,j}^{<t}(j)$ given the height $h_{i,L}(j)=\sup(e_{i,L}^{<t}(j))$ is $\Lambda_{h_{i,L}(j)}^p$. Putting all the above results together we have that given $t_i$ the random set $\mathcal{L}_i=\{(t_{i,L}(j),\Upsilon_{i,L}(j))\}_{j\geq 0}$ is conditionally independent of the set $\mathcal{R}_i$, and its law is a Poisson point-process with intensity measure

$$\mathbb{1}_{\{0<t<t_i\}}\,dt\,\mathbb{1}_{\{0<h<t\}}\frac{dh}{(1+h)^2}\Lambda_h^p.$$

The same conditional law of $\mathcal{R}_i$ follows from time reversibility of the law of $e_i^{<t}$.  □

Let us now consider the implications that the $p$-sampling of extinct individuals has in the asymptotic context. In Section 2, the genealogical point-process was defined from the contour process $\mathcal{C}_{\mathcal{T}_{t,n}}$, and its asymptotics was identified as the continuum genealogical point-process similarly defined from a Brownian excursion $\mathcal{B}_{t,1}$ conditioned to have local time 1 at level $t$. Now the $p$-sampled historical process is defined from a contour process $\mathcal{C}_{\mathcal{T}_{t,n}}$ whose local maxima are sampled independently with equal chance $p$. In terms of the (horizontal) $u$-coordinate of $\mathcal{C}_{\mathcal{T}_{t,n}}$ the $p$-sampled individuals form a random set of marks on $\mathbb{R}^+$. The fact that $\mathcal{C}_{\mathcal{T}}$ is an alternating sum of independent Exponential(rate 1) random variables implies that the random set formed by the local maxima of $\mathcal{C}_{\mathcal{T}}$ is a Poisson process of rate $1/2$ on $\mathbb{R}^+$, and the same still holds for the sets formed by the local maxima of each part of an excursion of $\mathcal{C}_{\mathcal{T}_{t,n}}$ below $t$. If we further sample these local maxima independently with chance $p$ we have a Poisson process of rate $p/2$ on $\mathbb{R}^+$. Now, for the asymptotics, the appropriate rescaling, as in Section 2, speeds up the time axis of $\mathcal{C}_{\mathcal{T}_{t,n}}$ by $n$. Hence if we consider $p_n$ such that $np_n \to p$ as $n \to \infty$, then asymptotically the $p_n$-sampling on $\mathcal{C}_{\mathcal{T}_{t,n}}$ will converge to a Poisson process of rate $p/2$. This prompts us to consider for the asymptotics of the $p$-historical point-process a process similarly defined from a conditioned Brownian excursion $\mathcal{B}_{t,1}$ sampled according to a Poisson(rate $p/2$) process along its (horizontal) $u$-coordinate.

REMARK 8. We are interested in obtaining an asymptotic point-process that has a.s. finitely many extinct individuals recorded. It is clear that thus the rate of sampling asymptotically has to satisfy $np_n \to p$ as $n \to \infty$.

We define a process derived from a conditioned Brownian excursion $\mathcal{B}_{t,1}$ in the same manner that $\Xi_{t,n}^p$ was derived from the contour process of the conditioned branching process $\mathcal{C}_{\mathcal{T}_{t,n}}$. Recall that $\mathcal{B}(u)$, $u \geq 0$, denotes a Brownian excursion, for a fixed $t > 0$; $\ell_t(u)$, $u \geq 0$, is its local time at level $t$ up to time $u$; $i_t(\ell)$, $\ell > 0$, is the inverse process of $\ell_t$. Also, $\mathcal{B}_{t,1}(u)$, $u \geq 1$ denotes the excursion $\mathcal{B}$ conditioned to have total local time at $t$ equal to 1, and $(\ell, e_\ell^{<t})$ denotes the set of excursions of $\mathcal{B}_{t,1}$ below level $t$ indexed by the local time $\ell_t$ at the time of their beginning.

Define the $p$-sampling on $\mathcal{B}_{t,1}$ to be a Poisson(rate $p/2$) process along the $u$-axis of $\mathcal{B}_{t,1}$. We indicate this by putting a star mark on the graph of $\mathcal{B}_{t,1}$ at the times of this Poisson process. Let $e_\ell^{<t}$ be an excursion of $\mathcal{B}_{t,1}$ below level $t$

$$e_\ell^{<t}(u) = \mathcal{B}_{t,1}(i_t(\ell^-) + u), \qquad u \in [0, i_t(\ell) - i_t(\ell^-)).$$

Recall that $a_\ell = \inf(e_\ell^{<t})$ is its lowest point occurring at $u_\ell = \arg\min(e^{<t}(u))$, and that $t_\ell = t - a_\ell$ denotes its distance from level $t$. For each $e_\ell^{<t}$ we define its left and right parts (relative to its lowest point) to be

$$e_{\ell,L}^{<t}(u) = \mathcal{B}_{t,1}(i_t(\ell^-) + u), \qquad u \in [0, u_\ell - i_t(\ell^-)),$$
$$e_{\ell,R}^{<t}(u) = \mathcal{B}_{t,1}(i_t(\ell) - u), \qquad u \in [0, i_t(\ell) - u_\ell).$$

Note that $e_{\ell,L}^{<t}$ runs forward to the lowest point of $e_\ell^{<t}$, whereas $e_{\ell,R}^{<t}$ runs backward in time to it. We shall also need their respective processes of infima

$$\varsigma_{\ell,L}(u) = \inf_{0 \leq v \leq u} e_{\ell,L}^{<t}(v), \qquad u \in [0, u_\ell - i_t(\ell^-)),$$

$$\varsigma_{\ell,R}(u) = \inf_{0 \leq v \leq u} e_{\ell,R}^{<t}(v), \qquad u \in [0, i_t(\ell) - u_\ell).$$

Figure 7 shows $e_{\ell,L}^{<t}$ and $e_{\ell,R}^{<t}$ with $\varsigma_{\ell,L}$ and $\varsigma_{\ell,R}$.

We define $a_{\ell,L}(j)$, $j \geq 0$, to be the successive levels of constancy of $\varsigma_{\ell,L}$, and we let $t_{\ell,L}(j) = t - a_{\ell,L}(j)$ be their distance to level $t$. For each level of constancy $a_{\ell,L}(j)$, let $e_{\ell,L}^{<t}(j)$ be the excursion of $e_{\ell,L}^{<t} - \varsigma_{\ell,L}$ that lies above the level $a_{\ell,L}(j)$. Let $h_{\ell,L}(j) = \sup(e_{\ell,L}^{<t}(j))$ be the height of this excursion. Note that a.s. all the $p$-sampled points on $\mathcal{B}_{t,1}$ lie on these excursions $e_{\ell,L}^{<t}(j)$. We define a tree $\Upsilon_{\ell,L}(j)$ induced by such a $p$-sampled excursion $e_{\ell,L}^{<t}(j)$ as the tree whose contour process is the linear interpolation of the sequence of the values of $e_{\ell,L}^{<t}(j)$ at the $p$-sampling times, alternating with the sequence of the minima of $e_{\ell,L}^{<t}(j)$ between the $p$-sampling times. An analogous definition leads to $a_{\ell,R}(j), j \geq 0$, $t_{\ell,R}(j)$, $j \geq 0$, $h_{\ell,R}(j)$ and $\Upsilon_{\ell,L}(R)$ from $e_{\ell,R}^{<t}(j)$ and $\varsigma_{\ell,R}(j)$.

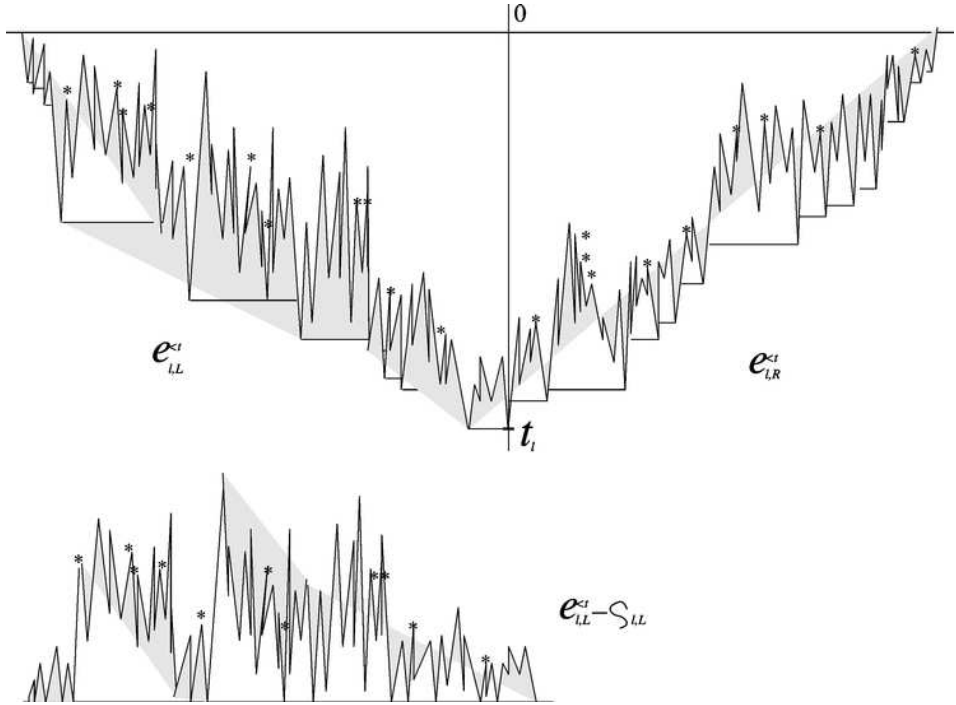

FIG. 7. (Top) An excursion $e_\ell^{<t}$ of $\mathcal{B}_{t,1}$ below $t$, its left $e_{\ell,L}^{<t}$ and right $e_{\ell,R}^{<t}$ parts, with their infimum processes; (bottom) the process $e_{\ell,L}^{<t} - \varsigma_{\ell,L}$.

REMARK 9. This definition of a tree from an excursion path sampled at given times has been explored for different sampling distributions in the literature (for some examples see [16], Section 6). Since for each $e_\ell^{<t}$ there are a.s. only finitely many $p$-sampled points the trees $\{\Upsilon_{\ell,L}(j)\}_j, \{\Upsilon_{\ell,R}(j)\}_j$ are a.s. in the space $\mathbf{T}$ of rooted planar trees with edge-lengths and finitely many leaves.

With each point $(\ell, t_\ell)$ of $\pi_{t,1}$ we now associate the sets

(16) $\quad \mathcal{L}_\ell = \{(t_{\ell,L}(j), \Upsilon_{\ell,L}(j))\}_{j \geq 0} \quad \text{and} \quad \mathcal{R}_\ell = \{(t_{\ell,R}(j), \Upsilon_{\ell,R}(j))\}_{j \geq 0}.$

We also define the first "right" set $\mathcal{R}_0$ and the last "left" set $\mathcal{L}_1$ from paths $e_{0,R}^{<t}$ of $\mathcal{B}_{t,1}$ before the first hitting time of $t$, and $e_{1,L}^{<t}$ of $\mathcal{B}_{t,1}$ after the last hitting time of $t$. For ease of notation we let $\mathcal{L}_0 = \mathcal{R}_1 = \varnothing$, $t_0 = t_1 = t$.

DEFINITION 10. The *$p$-sampled continuum historical point-process* $\xi_{t,1}^p$ is the random set

(17) $\qquad \xi_{t,1}^p = \{(\ell, t_\ell, \mathcal{L}_\ell, \mathcal{R}_\ell) : (\ell, t_\ell) \in \pi_{t,1},\, i_t(\ell^-) \neq i_t(\ell)\}.$

We next derive the law of the point-process $\xi_{t,1}^p$. For this we shall also need the law of the trees induced by the $p$-sampled excursions of $e^{<t} - \varsigma$. Let $\lambda^p$ denote the law on the space $\mathbf{T}$ induced by a $\mathcal{B}$ sampled at Poisson(rate $p$) points (in the sense of the bijection between sampled continuous functions and trees [3], same as the definition of $\Upsilon_{\ell,L}(j)$ from the $p$-sampled $e_{\ell,L}(j)$). Then, for any $h > 0$, let $\lambda_h^p$ denote the law induced by restricting $\lambda^p$ to the set of Brownian excursions $\mathcal{B}$ of height $h$.

To derive the law of $\xi_{t,1}^p$ we exploit in a more detailed manner the nice properties of Brownian excursions. We first use the result of Lemma 4, which gives us the law of the set $\{(\ell, t_\ell) : i_t(\ell^-) \neq i_t(\ell)\}$. Then conditional on this set we give the law of the sets $\mathcal{L}_\ell$ and $\mathcal{R}_\ell$. We show that $\{\mathcal{L}_\ell, \mathcal{R}_\ell\}_\ell$ are independent Poisson point-processes. The intensity measure of each such set is given by the following. First, choose $t_{\ell,L}(j)$, the distances below $t$ at which the $p$-sampled subtree excursions of $e_{\ell,L}^{<t} - \varsigma_{\ell,L}$ occur uniformly over $t_\ell$, the distance below $t$ of the lowest point of $e_\ell^{<t}$. Next, choose $h_{\ell,L}(j)$, the height for each such $p$-sampled excursion, according to the same law as that of the height of a $\mathcal{B}$ whose height is known to be less than $t_{\ell,L}(j)$. Finally, choose the law of the induced tree $\Upsilon_{\ell,L}(j)$ according to the law $\lambda_h^p$ described above.

LEMMA 7. *The random set $\xi_{t,1}^p$ is such that the following hold:*

(a) $\{(\ell, t_\ell) : i_t(\ell^-) \neq i_t(\ell)\}$ *is the Poisson point-process $\pi_{t,1}$ of Lemma* 4;

(b) *given* $\{(\ell, t_\ell) : i_t(\ell^-) \neq i_t(\ell)\}$ *the sets $\mathcal{L}_\ell$ and $\mathcal{R}_\ell$ are independent; and for each $\ell : i_t(\ell^-) \neq i_t(\ell)$, $\mathcal{L}_\ell$ and $\mathcal{R}_\ell$ are Poisson point-processes on $\mathbf{R}^+ \times \mathbf{T}$ with intensity measure*

$$\mathbb{1}_{\{0 < t < t_\ell\}} \, dt \, \mathbb{1}_{\{0 < h < t\}} \frac{dh}{h^2} \lambda_h^p. \tag{18}$$

PROOF. The proof proceeds in many of the same steps as the one for deriving the law of the $p$-sampled historical process $\Xi_{t,n}^p$. The notable difference is that we now have to resort to more sophisticated Markovian results on the decomposition of a Brownian path, such as the Williams decomposition of a Brownian excursion given its height, and the Pitman theorem on Bessel processes. In short, we consider the decomposition of the conditioned Brownian excursion $\mathcal{B}_{t,1}$ into its excursions below level $t$ provided by the Lemma 4. For each such excursion below $t$ given its lowest point at distance $t_\ell$ below $t$, the Williams decomposition gives us the independence and identity in law of its left and right parts, as well as the description of their laws in terms of a three-dimensional Bessel process. Furthermore, we can use Pitman's theorem that describes the law of the excursions of this Bessel process above the levels of constancy of its future infimum. After taking care of some conditioning issues, this finally gives us a simple description of these excursions above the levels of constancy as simply Brownian excursions conditioned on their maximal height.

The independence of the sets $\mathcal{L}_\ell$ over the index $\ell$ (the same holds for the sets $\mathcal{R}_\ell$) follows from the independence of the excursions of $\mathcal{B}_{t,1}$ below level $t$. This also holds (by the strong Markov property of $\mathcal{B}$) for the sets $\mathcal{R}_0$ and $\mathcal{L}_\infty$ defined from the parts of the path of $\mathcal{B}_{t,1}$ of its ascent to level $t$ and its descent from it. For each $e_\ell^{<t}$ excursion of $\mathcal{B}_{t,1}$ below level $t$, we let $e_\ell^+ = t - e_\ell^{<t}$. By Lemma 4, the conditional law of $e_\ell^+$ given $(\ell, t_\ell)$ is that of a Brownian excursion $\mathcal{B}$ conditioned on the value of its supremum $\mathcal{B} | \{\sup(\mathcal{B}) = t_\ell\}$. Let $\tau_{t_\ell} = \inf\{u > 0 : e_\ell^+(u) = t_\ell\}$; then by the Williams decomposition of a Brownian excursion $\mathcal{B}$ (e.g., [19], **1**, Section III.49), the law of $e_{\ell,L}^+ = t - e_{\ell,L}^{<t}$ is that of a Bess(3) (three-dimensional Bessel) process $\rho$ stopped the first time $\tau_{t_\ell}^\rho = \inf\{u > 0 : \rho(u) = t_\ell\}$ it hits $t_\ell$. By time reversibility of $\mathcal{B}$ the process

$$r_{\ell,L}(u) = t_\ell - e_{\ell,L}^+(\tau_{t_\ell} - u), \qquad u \in (0, \tau_{t_\ell}),$$

also has the law of the stopped Bess(3) process $\rho(u)$, $u \in (0, \tau_{t_\ell}^\rho)$. Let

$$j_{\ell,L}(u) = \inf_{u \leq v \leq \tau_{t_\ell}} r_{\ell,L}, \qquad u \in (0, \tau_{t_\ell}).$$

Then $\{t_\ell - t_{\ell,L}(j)\}_j$ are (in reversed index order) the successive levels of constancy of the process $j_{\ell,L}(u)$, $u \in (0, \tau_{t_\ell})$; $\{h_{\ell,L}(j)\}_j$ (in reversed index order) are the heights of the successive excursions from 0 of the process $r_{\ell,L}(u) - j_{\ell,L}(u)$, $u \in (0, \tau_{t_\ell})$, and $\{\Upsilon_{\ell,L}(j)\}_j$ (in reversed index order) are the trees induced by the $p$-sampled points on these excursions. To obtain the law of $j_{\ell,L}$ and $r_{\ell,L} - j_{\ell,L}$ consider the Bess(3) process $\rho(u), u \geq 0$, and its future infimum process $\jmath(u) = \inf_{v \geq u} \rho(v)$, $u \geq 0$. We note that the law of $j_{\ell,L}(u)$, $u \in (0, \tau_{t_\ell})$, is equivalent to that of $\jmath(u)$, $u \in (0, \tau_{t_\ell}^\rho)$, if $\jmath(\tau_{t_\ell}^\rho) = t_\ell$; in other words, if $\rho(u)$, $u \geq 0$, after it first reaches $t_\ell$ never returns to that height again. So,

$$(j_{\ell,L}, r_{\ell,L} - j_{\ell,L}) \stackrel{d}{=} (\jmath, \rho - \jmath) | \{\jmath(\tau_{t_\ell}^\rho) = t_\ell\} \qquad \text{for } u \in (0, \tau_{t_\ell}).$$

By Pitman's theorem, then by Levy's theorem (e.g., [18], VI, Sections 3 and 6)

$$(\jmath, \rho - \jmath) \stackrel{d}{=} (\zeta, \zeta - \beta) \stackrel{d}{=} (\bar{\ell}, |\bar{\beta}|),$$

where $\beta$ is a standard Brownian motion, $\zeta$ its supremum process; $|\bar{\beta}|$ is a reflected Brownian motion, $\bar{\ell}$ its local time at 0 (with the occupation time normalization). Thus, for $\bar{\tau}_{t_\ell} := \inf\{u \geq 0 : |\bar{\beta}|_u + \bar{\ell}_u = t_\ell\}$,

$$(j_{\ell,L}, r_{\ell,L} - j_{\ell,L}) \stackrel{d}{=} (\bar{\ell}, |\bar{\beta}|) | \{\bar{\ell}_{\bar{\tau}_{t_\ell}} = t_\ell\} \qquad \text{for } u \in (0, \tau_{t_\ell}).$$

The condition $\{\bar{\ell}_{\bar{\tau}_{t_\ell}} = t_\ell\}$ is equivalent to the condition $\{\bar{\ell}_{\bar{\tau}_{t_\ell}} = t_\ell, |\bar{\beta}|_{\bar{\tau}_{t_\ell}} = 0\}$ and $\{u < \bar{\tau}_{t_\ell} : \bar{\ell}_u < t_\ell, |\bar{\beta}|_u < t_\ell - \bar{\ell}_u\}$. Hence,

(19) $\quad (j_{\ell,L}, r_{\ell,L} - j_{\ell,L}) \stackrel{d}{=} (\bar{\ell}, |\bar{\beta}|) | \{\bar{\ell}_u < t_\ell, |\bar{\beta}|_u < t_\ell - \bar{\ell}_u; \bar{\ell}_{\bar{\tau}_{t_\ell}} = t_\ell, |\bar{\beta}|_{\bar{\tau}_{t_\ell}} = 0\}.$

Since $(\bar{\ell}, \sup(|\bar{\beta}|))$ is a Poisson point-process with intensity measure $d\bar{\ell}\, d\bar{h}/\bar{h}^2$, then using the independence property of a Poisson random measure on disjoint sets in (19), we obtain for $t = t_\ell - \bar{\ell}$ that $(t_\ell - j_{\ell,L}, \sup(r_{\ell,L} - j_{\ell,L}))$ is a Poisson point-process with intensity measure

$$\mathbb{1}_{(0<t<t_\ell)}\, dt\, \mathbb{1}_{(0<h<t)} \frac{dh}{h^2}.$$

Recall the relationship of the values $\{t_{\ell,L}(j), h_{\ell,L}(j), \Upsilon_{\ell,L}(j)\}_j$ of $\mathcal{L}_\ell$ with the processes $j_{\ell,L}$ and $r_{\ell,L} - j_{\ell,L}$. The above result thus implies that $\mathcal{L}_\ell$ is a Poisson point-process with intensity measure

$$\mathbb{1}_{(0<t<t_\ell)}\, dt\, \mathbb{1}_{(0<h<t)} \frac{dh}{h^2} \lambda_h^p,$$

where the last factor comes from the fact that $\Upsilon_{\ell,L}(j)$ is just the tree induced by the $p$-sampled excursion of $|\bar{\beta}|$ of height $h_{\ell,L}(j)$. $\square$

Our next goal is to show that the process $\xi_{t,1}^p$ whose law we have just obtained is indeed the asymptotic result of the processes $\Xi_{t,n}^p$ after appropriate rescaling. To do so, we first must show that the laws $\Lambda_h^{p_n}$ on the space of trees converge as $n \to \infty$ to the law $\lambda_h^p$ if $np_n \to p$. We need to consider more closely the trees $\Upsilon_{i,L}(j)$ and $\Upsilon_{\ell,L}(j)$ induced by the sampled excursions appearing in the historical point-processes above. In both cases we have an excursion, $\mathcal{C}_\mathcal{T}$ or $\mathcal{B}$, of a given height and with marks on it produced by a sampling process. Laws of the trees induced by sampled excursions of unrestricted height can be very simply and elegantly described (see [10] for the case of $\mathcal{B}$). However, for the trees from excursions of a given height that we need to consider here, the description is much messier. We shall give next a recursive description that applies equally to define an $\Upsilon_{l,L}(j)$ from $\mathcal{C}_\mathcal{T}$ of a given height, or to define $\Upsilon_{\ell,L}(j)$ from $\mathcal{B}$ of a given height. A similar recursive description of an infinite tree induced by an unsampled Brownian excursion is given by Abraham and Mazliak [2].

Define the "spine" of the tree to extend from the root of the tree to the point of maximal height in the excursion. An equivalent representation of the tree is one in which the subtrees of the trees on the left and on the right of the axis through the spine are attached to this spine, an example of which is shown in Figure 8. We obtain the branch levels at which these subtrees are attached, as well as parameters needed for the description of the subtrees as follows.

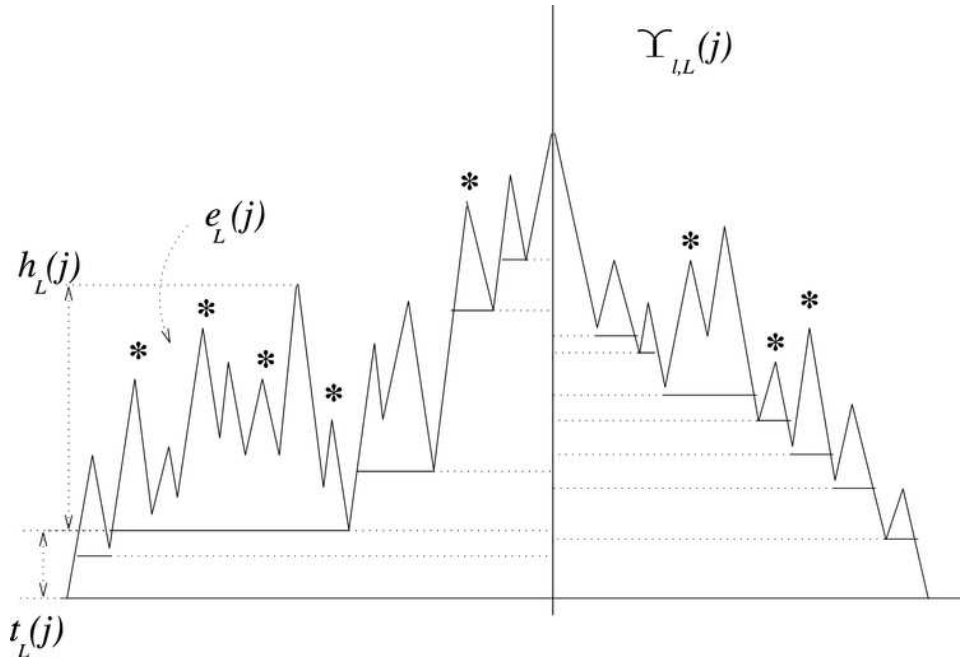

Fig. 8. *The "first" set in the recursive description consists of branch levels $\{t_L(j)\}_j$ at which subtrees induced by sampled excursions of $e_L - \varsigma_L$ are attached to the spine; and the heights $\{h_L(j)\}_j$ of these subtrees.*

We denote the excursion function defining this tree by $e(u)$, $u \geq 0$ (in other words $e = \mathcal{C}_\mathcal{T}$ or $e = \mathcal{B}$). Let $h$ be its given height, and $U_h = \arg\max\{e(u) : u > 0\}$ the time at which it is achieved. Then let $e_L(u)$, $u \in [0, U_h]$, be the left part of the excursion, and we also define its future infimum process $\varsigma_L(u) = \inf_{v \geq u} e(v)$, $u \in [0, U_h]$. Then the subtrees attaching on the left of the spine are defined by the process $e_L - \varsigma_L$ and the set of sampled marks. They are precisely the trees induced by the sampled excursions $e_L(j)$ of $e_L - \varsigma_L$ whose height is some $h_L(j)$. The levels at which they are attached to the spine are the levels of constancy $t_L(j)$ of $\varsigma_L$ at which the excursions of $e_L - \varsigma_L$ occur. Thus the set $\{(t_L(j), h_L(j))\}_{j \geq 0}$ is the "first" set in our recursive definition of trees. The "second" set is derived in the same manner from the sampled excursions $\{e_L(j)\}_j$ and so on. We define these sets analogously for the right part of $e$.

This recursive procedure is clearly very similar to our definition of the left and right sets $\mathcal{L}_i, \mathcal{R}_i$ for $e_i^{<t}$ and $\mathcal{L}_\ell, \mathcal{R}_\ell$ for $e_\ell^{<t}$ as defined earlier. The main difference is that the subtrees here are defined from excursions above the levels of constancy of the future infimum process for $e$, whereas earlier they were defined from excursions above the levels of constancy of the past infimum process for $e_i^{<t}$ and $e_\ell^{<t}$. However, time inversion and reflection invariance of the transition function of $e$ will allow us to easily derive the laws of the "first" set of points here from the results of Lemmas 6 and 7. In the next lemma we give a recursive description of the law of $\Lambda_h^{p_n}$ and $\lambda_h^p$, and we show that we do have the convergence of the $\Lambda_h^{p_n}$ (appropriately rescaled) to $\lambda_h^p$ if $np_n \to p$.

LEMMA 8. *The law $\Lambda_h^{p_n}$ of a tree induced by a $p_n$-sampled contour process $\mathcal{C}_\mathcal{T}$ of a given height $h$ is such that the first sets of points $\{t_L(j), h_L(j)\}_j$ and $\{t_R(j), h_R(j)\}_j$ are independent Poisson point-processes with intensity measure*

$$(20) \qquad \frac{1}{\sqrt{p_n}} \mathbb{1}_{(0 < \tau < h)} \, d\tau \, \mathbb{1}_{(0 < \kappa < h - \tau)} \frac{d\kappa}{(1+\kappa)^2} \frac{1+\tau}{\tau}.$$

*The law $\lambda_h^p$ of a tree induced by a $p$-sampled Brownian excursion $\mathcal{B}$ of a given height $h$ is such that the first sets of points $\{t_L(j), h_L(j)\}_j$ and $\{t_R(j), h_R(j)\}_j$ are independent Poisson point-processes with intensity measure*

$$(21) \qquad \frac{1}{\sqrt{p}} \mathbb{1}_{(0 < \tau < h)} \, d\tau \, \mathbb{1}_{(0 < \kappa < h - \tau)} \frac{d\kappa}{\kappa^2}.$$

*Let $n^{-1}\Lambda_h^{p_n}$ be the law of the tree induced by a rescaled $p_n$-sampled contour process $\mathcal{C}_\mathcal{T}$ by $n^{-1}$ in the vertical coordinate. Then for any $\{p_n \in (0,1)\}_{n \geq 1}$ such that $np_n \xrightarrow[n \to \infty]{} p$ we have $n^{-1}\Lambda_h^{p_n} \underset{n \to \infty}{\Longrightarrow} \lambda_h^p$.*

PROOF. The key for this proof is to observe the following. If $e(u)$, $u \geq 0$, is the $p_n$-sampled process $X_{\mathcal{T}}|\{\sup(\mathcal{C}_{\mathcal{T}}) = h\}$, then $e_L(u) = e(u)$, $u \in [0, U_h]$, has the law of a $p_n$-sampled $X_{\mathcal{T}}|\{\tau_h < \tau_0\}$, where $\tau_h, \tau_0$ are the first hitting times of $(h, +1), (0, -1)$, respectively, by $X_{\mathcal{T}}$. Then time reversibility and the reflection invariance of the transition function of $X_{\mathcal{T}}$ imply that $h - e_L(U_h - u)$, $u \in [0, U_h]$, has the same law as $e_L(u)$, $u \in [0, U_h]$. Now the levels of constancy of $\varsigma_L$, and the corresponding excursions $e_L - \varsigma_L$ above them, are equivalent to the levels of constancy and excursions of a set $\mathcal{L}_i$ considered in Lemma 6, thus giving a Poisson process of intensity measure as in (15). The factor $p^{-1/2}$ in the intensity measure (20) comes from the fact that here we only consider the excursions of $e_L - \varsigma_L$ that have at least one sampled mark in them. Namely, for the branching process $\mathcal{T}$, if $N_{\text{tot}}$ denotes the total population size of $\mathcal{T}$, then the generating function of $N_{\text{tot}}$ is $\mathbf{E}(x^{N_{\text{tot}}}) = 1 - (1 - x)^{1/2}$. Hence, the chance of at least one mark in the $p_n$-sampled point-process of $\mathcal{T}$ is $1 - \mathbf{E}((1 - p_n)^{N_{\text{tot}}}) = p_n^{1/2}$.

A similar argument applies when $e(u)$, $u \geq 0$, is the process $\mathcal{B}|\{\sup(\mathcal{B}) = h\}$ sampled at Poisson(rate $p/2$) times. Time reversibility and reflection invariance of the transition function of $\mathcal{B}$ allow us to identify that the law of the levels of constancy of $\varsigma_L$, and the corresponding excursions $e_L - \varsigma_L$ above them, are the same as those for a set $\mathcal{L}_\ell$ considered in Lemma 7, which we know form a Poisson process with intensity measure as in (18). The factor $p^{-1/2}$ in the intensity measure of (21) then comes from the rate of excursions with at least one sampled mark. Namely, a Poisson(rate $p/2$) process of marks on $\mathcal{B}$ along its time coordinate is in its local time coordinate a Poisson(rate $p^{1/2}$) process of marks (see [19], **2**, Section VI.50).

Now the law of the first set of the rescaled process with $n^{-1}\Lambda_h^{p_n}$ converges to the law of the first set of the process with the law $\lambda_h^p$. This follows from the fact that the former is a sequence of Poisson point-processes whose support set and intensity measure converge to those of the latter Poisson point-process. Since for Poisson random measures the convergence of finite-dimensional sets is sufficient to insure weak convergence of the whole process our claim follows for the first sets, and by recursion for the whole process. □

Finally, we can obtain the asymptotic result for the $p_n$-sampled historical point-processes. The rescaling of $\Xi_{t_n,n}^{p_n}$ is the same as that for $\Pi_{t,n}$. Both coordinates of $\Pi_{t,n}$ are rescaled by $n^{-1}$, so that the vertical coordinate of the sets $\mathcal{L}_i, \mathcal{R}_i$ is also rescaled by $n^{-1}$, and the sampling rate is rescaled by $n$. Hence the rescaled process is defined as

$$(22) \quad n^{-1}\Xi_{t_n,n}^{p_n} = \{(n^{-1}l_i, n^{-1}\tau_i, n^{-1}\mathcal{L}_i, n^{-1}\mathcal{R}_i) : (l_i, \tau_i, \mathcal{L}_i, \mathcal{R}_i) \in \Xi_{t_n,n}^{p_n}\}.$$

The asymptotic properties of the rescaled $p$-sampled historical process are now easily established from our earlier results.

THEOREM 9. *For any $\{t_n > 0\}_{n \geq 1}$ and $\{p_n \in (0,1)\}_{n \geq 1}$ such that $t_n/n \underset{n \to \infty}{\longrightarrow} t$ and $np_n \underset{n \to \infty}{\longrightarrow} p$ we have $n^{-1} \Xi_{t_n,n}^{p_n} \underset{n \to \infty}{\Longrightarrow} \xi_{t,1}^p$.*

PROOF. By Theorem 5 we already have that $n^{-1} \Pi_{t_n,n} \underset{n \to \infty}{\Longrightarrow} \pi_{t,1}$. Applying the rescaling to the results of Lemma 6 together with the result of Lemma 8 now implies that the support set and intensity measure of the Poisson point-process of each $\mathcal{L}_i$ after rescaling converges to those of the Poisson point-process $\mathcal{L}_\ell$ as given by Lemma 7. Then the convergence of the support set and intensity measure for the Poisson random measure $\Xi_{t_n,n}^{p_n}$ to those of $\xi_{t,1}^p$ implies the weak convergence of these processes. □

**Acknowledgments.** I thank D. Aldous for invaluable help and suggestions during the work on this paper. I also thank A. Winter for her thorough reading of this paper and helpful comments.

Department of Statistics
University of California
367 Evans Hall
Berkeley, California 94720
USA
e-mail: lea@stat.berkeley.edu
url: www.stat.berkeley.edu/users/lea